\documentclass[12pt]{amsart}
\usepackage{amssymb}
\usepackage{amsbsy}
\usepackage{amscd}
\usepackage[mathscr]{eucal}
\usepackage{verbatim}
\makeatletter
\def\cal{\mathcal}
\def\Bbb{\mathbb}
\def\frak{\mathfrak}

\newenvironment{pf*}[1]{\proof[#1]}{\endproof}
\makeatother
\newenvironment{aenume}{%
  \begin{enumerate}%
  }{\end{enumerate}}
\makeatletter
\@ifclasslater{amsart}{1999/11/24}{}{
\renewcommand*\subjclass[2][1991]{%
  \def\@subjclass{#2}%
  \@ifundefined{subjclassname@#1}{%
    \ClassWarning{\@classname}{Unknown edition (#1) of Mathematics
      Subject Classification; using '1991'.}%
  }{%
    \@xp\let\@xp\subjclassname\csname subjclassname@#1\endcsname
  }%
}
\renewcommand{\subjclassname}{%
  \textup{1991} Mathematics Subject Classification}
\@xp\let\csname subjclassname@1991\endcsname \subjclassname
\@namedef{subjclassname@2000}{%
  \textup{2000} Mathematics Subject Classification}
}
\makeatother
\newenvironment{NB}{
{\bf NB}. \footnotesize
}{}
  \renewenvironment{NB}{%
    \comment
    }{\endcomment}

\newtheorem{Theorem}[equation]{Theorem}
\newtheorem{Corollary}[equation]{Corollary}
\newtheorem{Lemma}[equation]{Lemma}
\newtheorem{Proposition}[equation]{Proposition}

\theoremstyle{definition}
\newtheorem{Definition}[equation]{Definition}

\newtheorem{Notation}[equation]{Notation}

\theoremstyle{remark}
\newtheorem{Remark}[equation]{Remark}

\numberwithin{equation}{section}
\newcommand{\thmref}[1]{Theorem~\ref{#1}}
\newcommand{\secref}[1]{\S\ref{#1}}
\newcommand{\lemref}[1]{Lemma~\ref{#1}}
\newcommand{\propref}[1]{Proposition~\ref{#1}}
\newcommand{\corref}[1]{Corollary~\ref{#1}}
\newcommand{\subsecref}[1]{\S\ref{#1}}

\newcommand{\defeq}{\overset{\operatorname{\scriptstyle def.}}{=}}
\newcommand{\C}{{\Bbb C}}
\newcommand{\Z}{{\Bbb Z}}
\newcommand{\Q}{{\Bbb Q}}
\newcommand{\R}{{\Bbb R}}
\newcommand{\proj}{{\Bbb P}}
\newcommand{\SU}{\operatorname{\rm SU}}
\newcommand{\GL}{\operatorname{GL}}

\newcommand{\algsl}{\operatorname{\frak{sl}}} % because \sl="slant"

\newcommand{\Spec}{\operatorname{Spec}\nolimits}

\newcommand{\Hom}{\operatorname{Hom}}
\newcommand{\Ext}{\operatorname{Ext}}

\newcommand{\rank}{\operatorname{rank}}

\newcommand{\pd}[2]{\frac{\partial#1}{\partial#2}}
\newcommand{\deriv}[1]{\frac{d}{d#1}}

\newcommand{\ve}{\varepsilon}

\newcommand{\linf}{\ell_\infty}
\newcommand{\shfO}{\mathcal O}
\newcommand{\bp}{{\widehat\proj}^2}
\newcommand{\bM}{{\widehat M}}

\newcommand{\Pic}{\operatorname{Pic}}

\newcommand{\Supp}{\operatorname{Supp}}
\newcommand{\ch}{\operatorname{ch}}
\newcommand{\Wedge}{{\textstyle \bigwedge}}
\newcommand{\Todd}{\operatorname{Todd}}
\newcommand{\Coh}{\operatorname{Coh}}

\newcommand{\Li}{\operatorname{Li}}
\newcommand{\Zin}{Z^{\text{\rm inst}}}
\newcommand{\bZin}{\widehat{Z}^{\text{\rm inst}}}
\newcommand{\bZ}{\widehat{Z}}

\newcommand{\Fin}{F^{\text{\rm inst}}}

\newcommand{\q}{\mathfrak q}
\newcommand{\bbeta}{\boldsymbol\beta}
\newcommand{\hT}{\widetilde T}

\begin{document}
\title[Instanton counting on blowup. II]
{Instanton counting on blowup. II.
\\ $K$-theoretic partition function}
\author{Hiraku Nakajima}
\address{Department of Mathematics, Kyoto University, Kyoto 606-8502,
Japan}
\email{nakajima@math.kyoto-u.ac.jp}
\thanks{The first author is supported by the Grant-in-aid
for Scientific Research (No.13640019, 15540023), JSPS}

\author{K\={o}ta Yoshioka}
\address{Department of Mathematics, Faculty of Science, Kobe University,
Kobe 657-8501, Japan}
\email{yoshioka@math.kobe-u.ac.jp}

\dedicatory{Dedicated to Vladimir Drinfeld on his fiftieth birthday}

\begin{abstract}
We study Nekrasov's deformed partition function
$Z(\ve_1,\ve_2,\vec{a};\q,\bbeta)$ of $5$-dimensional supersymmetric
Yang-Mills theory compactified on a circle. Mathematically it is the
generating function of the characters of the coordinate rings of the
moduli spaces of instantons on $\mathbb R^4$.
We show that it satisfies a system of functional equations, called
blowup equations, whose solution is unique.
As applications, we prove (a) $F(\ve_1,\ve_2,\vec{a};\q,\bbeta) =
\ve_1\ve_2 \log Z(\ve_1,\ve_2,\vec{a};\q,\bbeta)$ is regular at $\ve_1
= \ve_2 = 0$ (a part of Nekrasov's conjecture), and (b) the
genus $1$ parts, which are first several Taylor coefficients of
$F(\ve_1,\ve_2,\vec{a};\q,\bbeta)$, are written explicitly in terms of
$\tau = d^2 F(0,0,\vec{a};\q,\bbeta)/ da^2$ in rank $2$ case.
\end{abstract}

\maketitle

\section*{Introduction}

In Part I of this paper \cite{part1}, we studied Nekrasov's partition
function \cite{Nek} for $\mathcal N=2$ supersymmetric gauge theory in
$4$-dimension (see also \cite{lecture}). It is defined as the
generating function of the integral of the equivariant cohomology
class $1$ of the framed moduli space $M(r,n)$ of torsion free sheaves
on $\proj^2$ with rank $r$, $c_2 = n$:
\begin{equation*}
   \Zin(\ve_1,\ve_2,\vec{a};\q)
   = \sum_{n=0}^\infty \q^n \int_{M(r,n)} 1.
\end{equation*}
Here $(r+2)$-dimensional torus $\hT$ acts naturally on $M(r,n)$, and
$\ve_1$, $\ve_2$, $\vec{a} = (a_1,\dots,a_r)$ are generators of
$H^*_{\hT}(\mathrm{pt}) = S^*(\operatorname{Lie}\hT)$.
(More precisely, this is the instanton part of the partition function.
We multiply it with the {\it perturbative part}. See
\subsecref{subsec:pert}.)
It can be considered as series of equivariant Donaldson invariants for
$\R^4$, and there are close relation to the ordinary Donaldson
invariants, such as blowup formulas, wall-crossing formulas
\cite{lecture,GNY}.

In this part II, we study a similar partition function, in which we
replace the integration in the equivariant cohomology by one in
equivariant $K$-theory:
\begin{equation*}
   \Zin(\ve_1,\ve_2,\vec{a};\q,\bbeta)
   = \sum_n (\q\bbeta^{2r} e^{-r\bbeta(\ve_1+\ve_2)/2})^n 
   \sum_i (-1)^i \ch H^i(M(r,n),\shfO).
\end{equation*}
We consider $e^{\bbeta a_\alpha}$, $e^{\bbeta \ve_1}$, $e^{\bbeta
\ve_2}$ as characters of $\hT$ here. The formal parameter $\bbeta$ is
introduced so that the $K$-theoretic partition function converges to
the homological one when $\bbeta\to 0$.
It is called the partition function of the $5$-dimensional
supersymmetric gauge theory compactified on a circle in the physics
literature, where the radius of the circle is $\bbeta$.

Nekrasov \cite{Nek} conjectured that 
\(
      \Fin(\ve_1,\ve_2,\vec{a};\q)
       = \ve_1\ve_2\log \Zin(\ve_1,\ve_2,\vec{a};\q)
\)
is regular at $\ve_1, \ve_2 = 0$, and
\(
    \Fin_0(\vec{a};\q) \defeq \Fin(0,0,\vec{a};\q)
\)
is the instanton part of the Seiberg-Witten prepotential for $\mathcal
N=2$ supersymmetric gauge theory \cite{SW} with gauge group
$\SU(r)$. The Seiberg-Witten prepotential is defined by certain period
integrals of hyperelliptic curves, the so-called Seiberg-Witten
curves. (See \cite{lecture} for detail.)
This was a mathematically well formulated conjecture, which is
similar to the mirror symmetry.
The conjecture for the homological partition function was proved
affirmatively by the authors \cite{part1} and Nekrasov-Okounkov
\cite{NO} independently.
Nekrasov also conjectured the same statements for the above $K$-theoretic
version, and the technique in \cite{NO} can be applied to that case also.

We study the $K$-theoretic partition function through the approach
taken in \cite{part1}, which we briefly describe now: We consider
similar correlation functions given by generating functions of
characters of cohomology groups of Donaldson divisors $\mu(C)$ on the
moduli spaces $\bM(r,k,n)$ on the blowup plane $\bp$. As a simple
application of the Atiyah-Bott-Lefschetz formula, the correlation
functions can be expressed by $\Zin$. (See \eqref{eq:blow-up1}.)
On the other hand, through a geometric study of moduli spaces on the
blowup, we prove vanishing of certain cohomology groups. (See
\thmref{thm:key2}.) Combining these two results, we get a system of
functional equations, called {\it blowup equations}, satisfied by
$\Zin$. It determines the coefficients of $\q^n$ in $\Zin$
recursively. It also implies the regularity of
$\Fin(\ve_1,\ve_2,\vec{a};\q,\bbeta)$ and we get a differential
equation for $\Fin_0(\vec{a};\q,\bbeta)$, as limits of the blowup
equations. (See \eqref{eq:Kcontact2}.) We call it {\it a contact term
equation}, according to the name of their homological version. The
contact term equation also determines $\Fin_0(\vec{a};\q,\bbeta)$
recursively.

The homological version of the contact term equation was much studied
in the physics literature, and the Seiberg-Witten prepotential
satisfies the equations (see \cite{GM3,Mar} and the reference
therein). By the uniqueness of its solution, the instanton part of the
Seiberg-Witten prepotential is equal to $\Fin_0(\vec{a};\q)$. This
was our proof of Nekrasov's conjecture in \cite{part1}.
It is natural to hope that the same proof can be given for the
$K$-theoretic version. But we do not find our $K$-theoretic contact
term equation in the physics literature, and do not know how to prove
this assertion at this moment except for $r=2$ case.
\begin{NB}
Let us mention $r=2$ case. \verb+(^_^)+ (May 20) 
\end{NB}

Although we do not give the proof of Nekrasov's conjecture, we think
that it is worthwhile to pursue our approach by various reasons:
\begin{enumerate}
 \item The blowup equations determine not only
       $\Fin_0(\vec{a};\q,\bbeta)$, but also several higher
       coefficients of the expansion of
       $\Fin(\ve_1,\ve_2,\vec{a};\q,\bbeta)$ at $\ve_1=\ve_2=0$. (See
       \secref{sec:genus1}.) 
     \item The geometric study of moduli spaces on the blowup is
       probably useful for the study of the $K$-theoretic version of
       Donaldson invariants.
\end{enumerate}
Higher coefficients are identified with higher genus Gromov-Witten
invariants for certain noncompact Calabi-Yau $3$-folds (see
\cite{Nek}, \cite[\S7]{lecture}), and appear in the wall-crossing
formula for Donaldson invariants \cite{GNY}. Thus they are equally
important as $\Fin_0(\vec{a};\q,\bbeta)$.

For the ordinary Donaldson invariants, the vanishing of first several
blowup coefficients was well-known, and was proved by the dimension
counting argument. Our proof of the blowup equation for the homological
partition function was given by the same idea.
The proof for the $K$-theoretic partition function is very different,
and we use the Kawamata-Viewheg vanishing theorem, a result from
complex algebraic geometry. But we hope that a similar result holds
for the $K$-theoretic version of Donaldson invariants, whose existence
is still conjectural.

\subsection*{Acknowledgement}
The authors are grateful to the referee for helpful suggestions and comments.

\section{$K$-theoretic partition function}

In this section, we define the $K$-theoretic version of Nekrasov's
partition function. We follow \cite{part1,lecture} for which the
reader can find more detail and the references.

\subsection{Definition of the partition function}

Let $M(r,n)$ denote the framed moduli spaces of torsion free sheaves
$(E,\Phi)$ on $\proj^2$ with rank $r$ and $c_2 = n$. Let
$M^{\operatorname{reg}}_0(r,n)$ be the open subvariety consisting of
locally free sheaves. Let $M_0(r,n)$ be the Uhlenbeck (partial)
compactification of $M^{\operatorname{reg}}(r,n)$, i.e.,
\begin{equation*}
   M_0(r,n) = \bigsqcup_{n'=0}^n 
    M_0^{\operatorname{reg}}(r,n')\times S^{n-n'}\C^2.
\end{equation*}
We can endow this space with the structure of an affine algebraic
variety so that there is a projective morphism
\begin{equation*}
   \pi\colon M(r,n)\to M_0(r,n).
\end{equation*}
The corresponding map between closed points can be identified with
\begin{equation*}\label{eq:map_pi}
   (E,\Phi) \longmapsto
   ((E^{\vee\vee},\Phi), \operatorname{Supp}(E^{\vee\vee}/E))\in
   M_0^{\operatorname{reg}}(r,n')\times S^{n-n'}\C^2.
\end{equation*}
where $E^{\vee\vee}$ is the double dual of $E$ and
$\operatorname{Supp}(E^{\vee\vee}/E)$ is the support of
$E^{\vee\vee}/E$ counted with multiplicities.

Let $T$ be the maximal torus of $\GL_r(\C)$ consisting of diagonal
matrices and let $\hT = \C^*\times\C^* \times T$. We define an action
of $\hT$ on $M(r,n)$ as follows: For 
$(t_1,t_2)\in \C^*\times\C^*$, let $F_{t_1,t_2}$ be an automorphism of 
$\proj^2$ defined by
\[
    F_{t_1,t_2}([z_0: z_1 : z_2]) = [z_0: t_1 z_1 : t_2 z_2].
\]
For $\operatorname{diag}(e_1,\dots,e_r)\in T$ let $G_{e_1,\dots,e_r}$
denote the isomorphism of $\shfO_{\linf}^{\oplus r}$ given by
\[
    \shfO_{\linf}^{\oplus r}\ni (s_1,\dots, s_r) \longmapsto
     (e_1 s_1, \dots, e_r s_r).
\] 
Then for $(E,\Phi)\in M(r,n)$, we define
\begin{equation}\label{eq:action}
    (t_1,t_2,e_1,\dots,e_r)\cdot (E,\Phi)
    = \left((F_{t_1,t_2}^{-1})^* E, \Phi'\right),
\end{equation}
where $\Phi'$ is the composite of homomorphisms
\begin{equation*}
   (F_{t_1,t_2}^{-1})^* E|_{\linf} 
   \xrightarrow{(F_{t_1,t_2}^{-1})^*\Phi}
   (F_{t_1,t_2}^{-1})^* \shfO_{\linf}^{\oplus r}
   \longrightarrow \shfO_{\linf}^{\oplus r}
   \xrightarrow{G_{e_1,\dots, e_r}} \shfO_{\linf}^{\oplus r}.
\end{equation*}
Here the middle arrow is the homomorphism given by the action.

In a similar way, we have a $\hT$-action on $M_0(r,n)$. The map
$\pi\colon M(r,n)\to M_0(r,n)$ is equivariant.

\begin{Notation}\label{not:module}
We denote by $e_\alpha$ ($\alpha=1,\dots, r$) the one dimensional
$\hT$-module given by
\begin{equation*}
   \hT\ni (t_1,t_2, e_1, \dots, e_r) \mapsto e_\alpha.
\end{equation*}
Similarly, $t_1$, $t_2$ denote one-dimensional $\hT$-modules. Thus
the representation ring $R(\hT)$ is isomorphic to $\Z[t_1^\pm,
t_2^\pm, e_1^\pm, \dots, e_r^\pm]$, where $e_\alpha^{-1}$ is the dual
of $e_\alpha$.
\end{Notation}

We denote the coordinates of $\operatorname{Lie}(\hT)$
by $\ve_1,\ve_2,a_1,\dots,a_r$ corresponding to
$t_1,t_2,e_1,\dots,e_r$. In our previous paper \cite{part1}, these are
generators of the equivariant cohomology group $H^*_{\hT}(pt)$ of a
single point.
We relate two sets of variables as
\(
   t_1 = e^{\bbeta\ve_1}, t_2 = e^{\bbeta\ve_2},
   e_\alpha = e^{\bbeta a_\alpha},
\) 
where $\bbeta$ is a parameter. We will define the $K$-theory partition
function so that it converges to the homological partition function when
$\bbeta\to 0$.

We define the instanton part of the partition function by
\begin{equation}\label{eq:instantonpart}
\begin{gathered}
   \Zin(\ve_1,\ve_2,\vec{a};\q,\bbeta)
   \defeq \sum_n (\q\bbeta^{2r} e^{-r\bbeta(\ve_1+\ve_2)/2})^n 
   Z_n(\ve_1,\ve_2,\vec{a};\bbeta)\\
   Z_n(\ve_1,\ve_2,\vec{a};\bbeta) \defeq
   \sum_i (-1)^i \ch H^i(M(r,n),\shfO).
\end{gathered}
\end{equation}
Here the character $\ch$ is a formal sum of weight spaces, which are
finite dimensional as shown in \cite[\S4]{part1}. We also have
\begin{equation*}
   Z_n(\ve_1,\ve_2,\vec{a};\bbeta)
   = % \sum_n (\q\bbeta^{2r} e^{-r\bbeta(\ve_1+\ve_2)/2})^n
   \sum_i (-1)^i \ch H^0(M_0(r,n), R^i\pi_*\shfO).
\end{equation*}
We will see that the higher direct image sheaves $R^i\pi_*\shfO = 0$
for $i > 0$ (\lemref{lem:GR}). Therefore we have
\begin{equation*}
   Z_n(\ve_1,\ve_2,\vec{a};\bbeta)
   = % \sum_n (\q\bbeta^{2r} e^{-r\bbeta(\ve_1+\ve_2)/2})^n
   \ch H^0(M_0(r,n), \shfO).
\end{equation*}
This definition has an advantage that it involves only the Uhlenbeck
compactification $M_0(r,n)$.

\subsection{Other descriptions of the partition function}\label{subsec:other}

Let $K^{\hT}(M(r,n))$ denote the Grothendieck group of
$\hT$-equivariant coherent sheaves on $M(r,n)$ and similarly for
$K^{\hT}(M_0(r,n))$. These are modules over the representation ring
$R(\hT)$ of the torus $\hT$. As in \ref{not:module}, we identify it
with the Laurent polynomial ring
$\Z[t_1^\pm,t_2^\pm,e_1^\pm,\dots,e_r^\pm]$. Since $M(r,n)$ is
nonsingular, $K^{\hT}(M(r,n))$ is isomorphic to the Grothendieck group of
$\hT$-equivariant vector bundles. The proper equivariant morphism
$\pi$ induces a homomorphism 
\(
   \pi_* \colon K^{\hT}(M(r,n)) \to K^{\hT}(M_0(r,n))
\)
by taking the alternating sum of higher direct image sheaves $\sum_i
(-1)^i R^i \pi_*$.

The fixed points $M_0(r,n)^{\hT}$ consist of the single
point $n[0]\in S^n\C^2 \subset M_0(r,n)$. Let $\iota_0$ denote the
inclusion map of the fixed point.
By the localization theorem for the K-theory due to Thomason \cite{Th}
(a prototype was given in \cite{AS}), it is known that the
homomorphism $\iota_{0*}$ is an isomorphism after the localization:
\begin{equation*}
   \iota_{0*}\colon
   \mathcal R \cong
   K^{\hT}(M_0(r,n)^{\hT})\otimes_{R({\hT})}\mathcal R
   \xrightarrow{\cong}
   K^{\hT}(M_0(r,n))\otimes_{R({\hT})}\mathcal R,
\end{equation*}
where $\mathcal R = \Q(t_1,t_2,e_1,\dots,e_m)$ is the quotient field of
$R({\hT})$. We have $\iota_{0*}^{-1} = \ch$, which is a consequence of
a trivial identity $\ch\circ\iota_{0*} = \operatorname{id}$. We thus
get
\begin{equation*}
   Z_n(\ve_1,\ve_2,\vec{a};\bbeta)
   = % \sum_n (\q\bbeta^{2r} e^{-r\bbeta(\ve_1+\ve_2)/2})^n
   \left(\iota_{0*}\right)^{-1} \pi_*(\shfO_{M(r,n)}).
\end{equation*}
Here we denote the element in $K^{\hT}(M(r,n))$ corresponding to
$\shfO_{M(r,n)}$ by the same symbol for brevity. We hope that these
two meanings can be distinguished from the content.

The fixed points $M(r,n)^{\hT}$ consist of $(E,\Phi) =
(I_1,\Phi_1)\oplus\cdots\oplus (I_r,\Phi_r)$ such that
\begin{aenume}
\item $I_\alpha$ is an ideal sheaf of $0$-dimensional subscheme
$Z_\alpha$ contained in $\C^2 = \proj^2\setminus\linf$.
\item $\Phi_\alpha$ is an isomorphism from $(I_\alpha)_{\linf}$ to the $\alpha$th factor of $\shfO_{\linf}^{\oplus r}$.
\item $I_\alpha$ is fixed by the action of $\C^*\times\C^*$, coming
  from that on $\proj^2$.
\end{aenume}

We parametrize the fixed point set $M(r,n)^{\hT}$ by an $r$-tuple of
Young diagrams $\vec{Y} = (Y_1,\dots,Y_r)$ so that the ideal
$I_\alpha$ is spanned by monomials $x^i y^j$ placed at $(i-1,j-1)$
outside $Y_\alpha$. The constraint is that the total number of boxes
\(
    |\vec{Y}| \defeq \sum_\alpha |Y_\alpha|
\)
is equal to $n$.

Let $\iota$ denote the inclusion map of the fixed point set. We have
\begin{equation*}
   \iota_{*}\colon
   \bigoplus_{\vec{Y}} \mathcal R \cong
   K^{\hT}(M(r,n)^{\hT})\otimes_{R({\hT})}\mathcal R
   \xrightarrow{\cong}
   K^{\hT}(M(r,n))\otimes_{R({\hT})}\mathcal R.
\end{equation*}
We then have
\begin{equation*}
   Z_n(\ve_1,\ve_2,\vec{a};\bbeta)
   = % \sum_n (\q\bbeta^{2r} e^{-r\bbeta(\ve_1+\ve_2)/2})^n
   \sum_{\vec{Y}}
   \left(\iota_{*}\right)^{-1}(\shfO_{M(r,n)}).
\end{equation*}

As $M(r,n)$ is nonsingular, $(\iota_*)^{-1}$ can be explicitly given
by Atiyah-Bott Lefschetz fixed points formula:
\begin{equation*}
   \iota_*^{-1}(\bullet) =
   \bigoplus_{\vec{Y}} \frac{\iota_{\vec{Y}}^*(\bullet)}
   {\Wedge_{-1}T^*_{\vec{Y}}M(r,n)},
\end{equation*}
where $T^*_{\vec{Y}}M(r,n)$ is the cotangent bundle of $M(r,n)$ at a
fixed point of $\vec{Y}$ considered as a ${\hT}$-module, $\Wedge_{-1}$ is
the alternating sum of exterior powers, and $\iota_{\vec{Y}}^*$ is the
pull-back homomorphism with respect to the inclusion
$\iota_{\vec{Y}}\colon \{ \vec{Y}\}\to M(r,n)$. Here the pull-back
homomorphism is defined via the isomorphism of $K^{\hT}(M(r,n))$ and the
Grothendieck group of ${\hT}$-equivariant {\it locally free sheaves\/}.

In order to express $\Wedge_{-1} T^*_{\vec{Y}}(M(r,n))$, we need some
notation for Young diagrams.
Let $Y = (\lambda_1\ge \lambda_2 \ge \cdots)$ be a Young diagram,
where $\lambda_i$ is the length of the $i$th column.
Let $Y' = (\lambda'_1\ge \lambda_2' \ge \dots)$ be the transpose of
$Y$. Thus $\lambda'_j$ is the length of the $j$th row of $Y$.
Let $l(Y)$ denote the number of columns of $Y$, i.e., $l(Y) =
\lambda'_1$.
Let
\begin{equation*}
a_Y(i,j) = \lambda_i - j, \qquad l_Y(i,j) = \lambda'_j - i.
\end{equation*}
% \begin{alignat*}{2}
% & a_Y(i,j) = \lambda_i - j, & \qquad & a'(i,j) = j - 1 \\
% & l_Y(i,j) = \lambda'_j - i, &\qquad & l'(i,j) = i - 1.
% \end{alignat*}
Here we set $\lambda_i = 0$ when $i > l(Y)$. Similarly $\lambda'_j =
0$ when $j > l(Y')$.
When the square $s = (i,j)$ lies in $Y$, these are called {\it
arm-length}, {\it leg-length}, respectively in the literature. But our
formula below involves these also for squares outside $Y$. So these
take negative values in general.

By \cite[Theorem~2.11]{part1} we have
\begin{equation}\label{eq:Zsum}
    \Zin(\ve_1,\ve_2,\vec{a};\q,\bbeta)
   = \sum_{\vec{Y}}
     \frac{(\q \bbeta^{2r} e^{-r\bbeta(\ve_1+\ve_2)/2})^{|\vec{Y}|}} 
     {\Wedge_{-1}T_{\vec{Y}}^*M(r,n)}
   = \sum_{\vec{Y}} 
   \frac{(\q \bbeta^{2r} e^{-r\bbeta(\ve_1+\ve_2)/2})^{|\vec{Y}|}}
    {\displaystyle\prod_{\alpha,\beta}
    n^{\vec{Y}}_{\alpha,\beta}(\ve_1,\ve_2,\vec{a};\bbeta)},
\end{equation}
where
\begin{multline*}
n^{\vec{Y}}_{\alpha,\beta}(\ve_1,\ve_2,\vec{a};\bbeta)
   = \prod_{s \in Y_\alpha}
      \left(1-e^{-\bbeta(-l_{Y_\beta}(s)\ve_1 
        + (a_{Y_\alpha}(s)+1)\ve_2 + a_\beta - a_\alpha)}\right)
\\
    \times\prod_{t\in Y_\beta} 
     \left(1-e^{-\bbeta((l_{Y_\alpha}(t)+1)\ve_1
       -a_{Y_\beta}(t)\ve_2 + a_\beta  - a_\alpha)}\right)
.
\end{multline*}

\begin{Remark}
\begin{NB}
Kota does not take the convention . I need
to check the consistency carefully.
\end{NB}
Contrary to the convention in \cite[Remark~4.4]{part1}, we put the
$\hT$-module structure on the coordinate ring (and the cohomology
groups) by $F_{g^{-1}}^*$, where $F_g \colon M(r,n)\to M(r,n)$ is the
isomorphism given by an element $g\in \hT$.
\end{Remark}

This combinatorial expression was the original definition of the
instanton part of the partition function due to Nekrasov \cite{Nek},
except we put the additional factor $e^{-r\bbeta(\ve_1+\ve_2)/2}$.
This factor is the half of the canonical bundle of $M(r,n)$ (see
\lemref{lem:canonical}), hence $Z_n$ is the index of the Dirac
operator, rather than the Dolbeault operator. Also, the factor makes
the symmetry of the partition function nicer, as we see in
\lemref{lem:blow-up}.

Moreover, it is clear that this $K$-theoretic partition function
converges to the homological partition function studied in
\cite{part1,lecture} as $\bbeta\to 0$.

\section{Blowup equation and a main result}

\subsection{Correlation functions on blowup}

Let $\bp$ be the blowup of $\proj^2$ at $[1:0:0]$. Let $p\colon
\bp\to\proj^2$ denote the projection. Let $C$ be the exceptional
divisor, $\shfO(C)$ the corresponding line bundle, and $\shfO(mC)$ its
$m$th power.

Let $\bM(r,k,\widehat{n})$ be the framed moduli space of torsion free
sheaves $(E,\Phi)$ on $\bp$ with rank $r$, $\langle c_1(E),[C]\rangle
= -k$ and $\langle \Delta(E), [\bp]\rangle = \widehat{n}$ where
$\Delta(E) = c_2(E) - \frac{r-1}{2r} c_1(E)^2$. This is also
nonsingular of dimension $2\widehat{n}r$. (Remark that $\widehat{n}$
may not be integer in general.) By tensoring a line bundle if
necessary, we assume $0\le k<r$ hereafter.

By \cite[Theorem~3.3]{lecture} there is a projective morphism defined by
\begin{alignat*}{2}
   \widehat\pi\colon& \bM(r,k,\widehat{n}) &\quad \to\quad & M_0(r,n)
\\
   & (E,\Phi) & \quad \mapsto\quad & \left(((p_* E)^{\vee\vee}, \Phi),
   \Supp(p_*E^{\vee\vee}/p_*E) + \Supp(R^1p_* E)\right),
\end{alignat*}
where $n = \widehat{n}-\frac{k(r-k)}{2r}$.

Let $(\mathcal E,\Phi)$ be the universal family on
$\bp\times\bM(r,k,\widehat{n})$. It has a natural $\hT$-structure from
the construction. As $C$ is $\C^*\times\C^*$-invariant, its
fundamental class $[C]$ defines a class in the equivariant homology group.
We define an equivariant $\Q$-divisor $\mu(C)$ on
$\bM(r,k,\widehat{n})$ by
\begin{equation*}
   \mu(C) \defeq 
   p_{2*}\left(
   \Delta(\mathcal E)\cap
%   \left(c_2(\mathcal E) - \frac{r-1}{2r} c_1(\mathcal E)^2\right)\cap
   \left[C\times\bM(r,k,\widehat{n})\right]\right)
   \in A_{2\widehat{n}r-1}^{\hT}(\bM(r,k,\widehat{n}))_\Q,
\end{equation*}
where $p_2\colon\bp\times\bM(r,k,\widehat{n})\to\bM(r,k,\widehat{n})$
is the projection to the second factor. As customary, we denote the
above as the `slant product'
\( 
   \Delta(\mathcal E)/[C].
%   \left(c_2(\mathcal E) - \frac{r-1}{2r} c_1(\mathcal E)^2\right)/[C].
\)
In the ordinary nonequivariant situation, it is known that $\mu(C)$
comes from the determinant line bundle
\begin{equation*}
   \det\left(R^\bullet p_{2*}\left( \mathcal E_{|C}\otimes
       p_1^*\shfO_C(-1)\right)\right).
\end{equation*}
(See \cite[\S1]{Li} for more detail).
This construction works in our equivariant setting, so $\mu(C)$ comes
from an equivariant $\Q$-line bundle. We will identify $\mu(C)$ as the
latter element hereafter, but this does not make any trouble as
$\Pic^{\hT}(\bM(r,k,\widehat{n})) \cong
A^{\hT}_{2\widehat{n}r-1}(\bM(r,k,\widehat{n}))$ (see
\cite[Lemma~1.3]{Y} or \cite[Theorem~1]{EG}).

\begin{NB}
  When $c_1 = 0$, $-\mu(C)$ is the first Chern class of the
  determinant line bundle. If $c_1\neq 0$, we need a modification.
\end{NB}

We want to define the $K$-theoretic direct image 
$\widehat{\pi}_*({\cal O}(d\mu(C)))$ 
of the $\hT$-equivariant ${\Bbb Q}$-line bundle $\mu(C)$.
For this purpose, we compute
$\det{\cal E}$ for the universal family
$({\cal E},\Phi)$ on $\widehat{\Bbb P}^2 \times \bM(r,k,\widehat{n})$.  
Since $h^1(\widehat{\Bbb P}^2)=0$,
$\det{\cal E} \cong
L \boxtimes {\cal O}_{\widehat{\Bbb P}^2}(kC)$
for a $\hT$-line bundle $L$ on $\bM(r,k,\widehat{n})$.
\begin{NB}
Note that $\shfO(C)$ has a natural $\hT$-structure so that
$0 \to \shfO(-C) \to \shfO \to \shfO_C \to 0$ is equivariant.
\end{NB}
Since ${\cal E}_{|\linf \times \bM(r,k,\widehat{n})}
\cong \bigoplus_{\alpha=1}^r 
{\cal O}_{\linf \times \bM(r,k,\widehat{n})}e_\alpha$ as
$\hT$-sheaves, we have
$L \cong {\cal O}_{\bM(r,k,\widehat{n})}\prod_\alpha e_\alpha$.
Hence $c_1(\mathcal E) = \sum a_\alpha + kC$, therefore we have
\begin{equation*}
\begin{split}
\mu(C) & =
\left(c_2({\cal E})-\frac{r-1}{2r}\left(\sum_\alpha a_\alpha+kC \right)^2 
\right)
/[C]\\
& =
c_2({\cal E})/[C]+\lambda,
\end{split}
\end{equation*}
where $\lambda=\frac{r-1}{2r}
\left(2k\sum_\alpha a_\alpha+k^2(\ve_1+\ve_2)\right) \in H^2(B\hT,{\Bbb Q})$.
Here we have used $1/[C] = 0$, $C/[C] = -1$, $C^2/[C] =
-(\ve_1+\ve_2)$ (cf.\ \cite[Proof of Lemma~5.8]{part1}).

Now we define the $K$-theoretic direct image as follows:
\begin{Definition}
\begin{equation*}
\widehat{\pi}_*({\cal O}(d \mu(C))) \defeq
\widehat{\pi}_*({\cal O}(d c_2({\cal E})/[C])) 
\otimes {\cal O}_{M_0(r,n)}(d \lambda) 
\in K^{\hT}(M_0(r,n)) \otimes_{R(\hT)} {\cal R}',
\end{equation*}
where ${\cal R}'={\Bbb Q}[t_1^{1/r},t_2^{1/r},e_1^{1/r},...,e_r^{1/r}]$.
\end{Definition} 

Note that $\widehat{\pi}_*({\cal O}(d \mu(C)))$ is in
$K^{\hT}(M_0(r,n))$ if $k$ is $0$ or $r$, or $d = 0$.

We now define {\it correlation functions\/} on blowup:
\begin{equation*}
  \bZin_{k,d}(\ve_1,\ve_2,\vec{a};\q,\bbeta) \defeq
  \sum_{\widehat{n}} (\q \bbeta^{2r} e^{-r\bbeta(\ve_1+\ve_2)/2}
  )^{\widehat{n}} 
  \left(\iota_{0*}\right)^{-1}\left( \widehat\pi_* ({\cal O}(d
  \mu(C))) \right).
\end{equation*}

\subsection{Correlation functions via the partition function}
As in \subsecref{subsec:other} we can use Atiyah-Bott Lefschetz fixed
points formula to get another description of $\bZin_{k,d}$. The
necessary computation for the weights of tangent spaces and divisor
$\mu(C)$ at fixed points was already done in homological version
\cite[\S3, \S6]{part1}. So we only state the answer:
\begin{multline}\label{eq:blow-up1}
%\begin{equation}\label{eq:blow-up1}
%\begin{split}
 %&
  \bZin_{k,d}(\ve_1,\ve_2,\vec{a};\q,\bbeta) =
  \sum_{\{\vec{k}\}=-{k}/{r}}
   \frac{(e^{\bbeta(\ve_1+ \ve_2) (d - r/2)} \q
   \bbeta^{2r})^{(\vec{k},\vec{k})/2}e^{\bbeta(\vec{k},\vec{a})d}}
   {\prod_{\vec\alpha\in\Delta} l^{\vec{k}}_{\vec{\alpha}}(\ve_1,\ve_2,\vec{a})} \times\\
% &
   \Zin(\ve_1,\ve_2-\ve_1,\vec{a}+\ve_1\vec{k};
    e^{\bbeta\ve_1 (d-r/2)} \q,\bbeta)
   \Zin(\ve_1-\ve_2,\ve_2,\vec{a}+\ve_2\vec{k};
    e^{\bbeta\ve_2 (d-r/2)}\q,\bbeta).
%\end{split}
%\end{equation}
\end{multline}
We need explanations of several notations. The vector $\vec{a}$ is
considered as an element of the Cartan subalgebra $\mathfrak h$ of
$\algsl_r$ by imposing the condition $\sum_\alpha a_\alpha = 0$. Then
$\Delta$ is the set of roots of $\algsl_r$. The vector $\vec{k}$ runs
over the coweight lattice
\(
    P = \{ \vec{k} = (k_1,\dots,k_r)\in\Q^r \mid \sum k_\alpha = 0,
  \exists k\in\Z\; \forall \alpha\; k_\alpha \equiv - {k}/{r} \mod \Z
%  \text{there exists $k\in \Z$ such that $k_\alpha \equiv -\frac{k}r
%    \mod \Z$ for all $\alpha$}
  \}
\)
with the constraint
\(
   \{ \vec{k} \} \defeq k_\alpha\pmod\Z = - {k}/{r}.
\) 
Finally  we set
\begin{equation}\label{eq:l}
   l^{\vec{k}}_{\vec{\alpha}}(\ve_1,\ve_2,\vec{a}) = 
   \begin{cases}
     {\displaystyle
     \prod_{\substack{i,j\ge 0\\i+j \le -\langle\vec{k},
     \vec{\alpha}\rangle-1}}}
          (1-e^{\bbeta(i\ve_1 +j\ve_2 - \langle \vec{a}, \vec{\alpha}\rangle)})
       & \text{if $\langle \vec{k}, \vec{\alpha}\rangle < 0$}, \\
     {\displaystyle
     \prod_{\substack{i,j\ge 0\\i+j\le \langle \vec{k},
     \vec{\alpha}\rangle-2}}}
          \left(1-e^{\bbeta(-(i+1)\ve_1 - (j+1)\ve_2 - 
         \langle \vec{a},\vec{\alpha}\rangle)}\right)
       & \text{if $\langle\vec{k}, \vec{\alpha}\rangle > 1$}, \\
     1 & \text{otherwise}
   \end{cases}
\end{equation}
for a root $\vec{\alpha}\in\Delta$.

\begin{NB}
  The above calculation of weights is consistent with \cite{part1}.
  If $\{ w_i \}$ are weights of $E_{|(\mathrm{fixed point})}$, the Euler
  class is $\prod w_i$, and $\Wedge_{-1} E^*$ is $\prod (1 - e^{-w_i})$.
\end{NB}

\subsection{Main results}

We can now state our main results:
\begin{Theorem}\label{thm:key2}
\textup{(1)(}$d=0$ case\textup) 
\begin{equation}\label{eq:blow-up2}
\bZin_{k,0}(\ve_1,\ve_2,\vec{a};\q,\bbeta)
=
 (\q \bbeta^{2r} e^{-r\bbeta(\ve_1+\ve_2)/2})^{\frac{k(r-k)}{2r}}
\Zin(\ve_1,\ve_2,\vec{a};\q,\bbeta).
\end{equation}

\textup{(2)(}$0<d<r$ case\textup)
\begin{equation}\label{eq:blow-up3}
   \bZin_{k,d}(\ve_1,\ve_2,\vec{a};\q,\bbeta) = 
 \begin{cases}
   \Zin(\ve_1,\ve_2,\vec{a};\q,\bbeta) & \text{for $k = 0$,}\\
   0 & \text{for $0 < k < r$.}
 \end{cases}  
\end{equation}

\textup{(3)(}$d=r$ case\textup)
\begin{equation}
\bZin_{k,r}(\ve_1,\ve_2,\vec{a};\q,\bbeta)=
(-1)^{k(r-k)}(t_1 t_2)^{{k(r-k)}/{2}}
 (\q \bbeta^{2r} e^{-r\bbeta(\ve_1+\ve_2)/2})^{\frac{k(r-k)}{2r}}
\Zin(\ve_1,\ve_2,\vec{a};\q,\bbeta).
\end{equation}
\end{Theorem}

\begin{NB}
I put the factor coming from the difference of the instanton numbers.
(May 20)
\end{NB}

For the homological partition function, similar formulas were obtained
(\cite[6.12]{part1}, \cite[\S5.1]{lecture}). The proof was the same as
that of first several coefficients of the blowup formula for Donaldson
invariants. It essentially follows from the dimension counting
argument. We hope that the above formulas can be also considered as
blowup formulas in low degree for the $K$-theoretic equivariant
Donaldson invariants, whose existence is still conjectural.

Combining \thmref{thm:key2} with \eqref{eq:blow-up1}, we get
functional equations satisfied by $\Zin$. These equations, which we
call {\it blowup equations\/}, are powerful, and we have several
consequences as we will see in later sections. We give the first
application now:

\begin{Corollary}\label{cor:recursive}
Let $Z_n \equiv Z_n(\ve_1,\ve_2,\vec{a};\bbeta)$ be the coefficient of
$\q^n$ in $\Zin(\ve_1,\ve_2,\vec{a};\q,\bbeta)$ as in
\eqref{eq:instantonpart}. Then $Z_n$ is determined from $Z_0 = 1$
inductively by the blowup equations from \thmref{thm:key2} with
$k=0$, $d=0,1,2$. 
\end{Corollary}

The proof is contained in that of \thmref{thm:regular}.

\section{Vanishing theorem}

We will prove \thmref{thm:key2} in this section. The crucial result
here is the vanishing theorem of higher direct image sheaves (see
\propref{prop:vanish}). A reader who has interests only in
applications of equations in \thmref{thm:key2} can safely skip this
section.

\subsection{}

%For an $\hT$-equivariant ${\Bbb Q}$-divisor $D$ on
%${\Bbb P}^2$ or $\widehat{\Bbb P}^2$,
%$L_D$ denotes the ${\Bbb Q}$-line bundle such that
%$c_1(L_D)=\mu(D)$. 

%Let $\nu:M_0(r,n)_{nor} \to M_0(r,n)$ be the normalization 
%of $M_0(r,n)$. Since $M_0(r,n)$ is an affine variety,
%$M_0(r,n)_{nor}$ is an affine variety with the structure sheaf
%$\pi_*({\cal O}_{M(r,n)})$.
%Obviously $M_0(r,n)_{nor}$.
%has a compatible $\hT$-action with the normalization map.
%Moreover
%$\pi$ and 
%$\widehat{\pi}$ factor through $M_0(r,n)_{nor}$,
%and they are $\hT$-equivariant.
%\begin{NB}
%Uniqueness of the normalization.
%\end{NB}
%In particular
%$\widehat{\pi}_*({\cal O}_{\bM(r,k,n)})\cong
%\nu_*({\cal O}_{M_0(r,n)_{nor}}) \cong \pi_*({\cal O}_{M(r,n)})$
%as $\hT$-equivariant coherent sheaves.
%Since $\nu$ is bijective,
%replacing $M_0(r,n)$ by $M_0(r,n)_{nor}$, we may assume that 
%$M_0(r,n)$ is normal.
%\begin{NB}
%$M_0(r,n)$ will be normal.
%\end{NB}

It is known that $M_0(r,n)$ is a normal variety by \cite{CB}.
Since $\pi:M(r,n) \to M_0(r,n)$ is birational,
$\pi_*({\cal O}_{M(r,n)})={\cal O}_{M_0(r,n)}$.
Since $M(r,n)$ is a holomorphic symplectic manifold,
$K_{M(r,n)} \cong {\cal O}_{M(r,n)}$ as sheaves
(see Lemma~\ref{lem:canonical} for a different proof).
By the Grauert-Riemenschneider vanishing theorem,
the higher direct image sheaves vanish, i.e.,
$R^i \pi_*({\cal O}_{M(r,n)})=0$ for $i > 0$. Hence we have the
following.
\begin{Lemma}\label{lem:GR}
We have the following equality in the equivariant
$K$-group $K^{\hT}(M_0(r,n))$:
\begin{equation*}
\pi_{*}({\cal O}_{M(r,n)})={\cal O}_{M_0(r,n)}.
\end{equation*}
\end{Lemma}

\thmref{thm:key2} is equivalent to the following:
\begin{Proposition}\label{prop:key}
We have the following equalities 
in $K^{\hT}(M_0(r,n)) \otimes_{R(\hT)} {\cal R}$:

\textup{(1)}
\begin{equation*}
\widehat{\pi}_* 
({\cal O}_{\bM(r,k,\widehat{n})})={\cal O}_{M_0(r,n)}.
\end{equation*}

\textup{(2)}
If $0<d<r$, then
\begin{equation*}
\widehat{\pi}_* 
({\cal O}_{\bM(r,k,\widehat{n})}(d\mu(C)))=
\begin{cases}
{\cal O}_{M_0(r,n)},&\;k=0\\
0,& \;0<k<r.
\end{cases}
\end{equation*}

\textup{(3)}
\begin{equation*}
\widehat{\pi}_* 
({\cal O}_{\bM(r,k,\widehat{n})}(r\mu(C)))=
(-1)^{k(r-k)}(t_1 t_2)^{{k(r-k)}/{2}}
{\cal O}_{M_0(r,n)}.
\end{equation*}
\end{Proposition}

The equation for the case where $0 \le d <r$ is a consequence of the
following two propositions.
\begin{Proposition}\label{prop:image}
In the category of the equivariant coherent sheaves
$\Coh_{\hT}(M_0(r,{n}))$, we have the following
equalities: 

\textup{(1)}
\begin{equation*}
\widehat{\pi}_* 
({\cal O}_{\bM(r,k,\widehat{n})})=\pi_*({\cal O}_{M(r,n)}).
\end{equation*}

\textup{(2)}
If $d \geq 1$, then 
\begin{equation*}
\widehat{\pi}_* 
({\cal O}_{\bM(r,k,\widehat{n})}(d \mu(C)))=
\begin{cases}
 {\cal O}_{M_0(r,n)}, & k=0\\
 0, & 0<k<r.
\end{cases}
\end{equation*}
\end{Proposition}

\begin{Proposition}\label{prop:vanish}
\begin{equation*}
R^i\widehat{\pi}_* 
({\cal O}_{\bM(r,k,\widehat{n})}(d \mu(C)))=0
\end{equation*}
for $i>0$ and $d < r$.
\end{Proposition}

Proof of Proposition~\ref{prop:image}:
We start with the proof of (1). We note that 
$\widehat{\pi}$ is a Grassmannian bundle over
$M_0^{\mathrm{reg}}(r,n)$ with a fiber $Gr(r,k)$.
In fact, the inverse image of $(E,\Phi)\in M_0^{\mathrm{reg}}(r,n)$
under $\widehat{\pi}$ consists of sheaves $E'$, which fit in an exact
sequence
\[
   0 \to E'(-C) \to p^*E \to \shfO_C^{\oplus r-k} \to 0.
\]
\begin{NB}
Please check this is correct:
   $c_1(E') = c_1(E'(-C)) + r[C] = (k-r)[C] + r[C] = k[C]$.
(May 20)
\end{NB}
Thus ${\cal O}_{M_0(r,n)} \to 
\widehat{\pi}_*({\cal O}_{\bM(r,k,\widehat{n})})$
is an isomorphism on $M_0^{\mathrm{reg}}(r,n)$.
Hence $\Spec(\widehat{\pi}_*({\cal O}_{\bM(r,k,\widehat{n})}))
\to M_0(r,n)$ is finite and birational.
Since $M_0(r,n)$ is normal, the assertion holds.
\begin{NB}
Only the connectivity of general fibers of $\widehat{\pi}$ (*) is
essential:
We take the Stein factorization:
$\bM(r,k,\widehat{n}) \to Y \to M_0(r,n)$, where
$Y:=\Spec(\widehat{\pi}_*({\cal O}_{\bM(r,k,\widehat{n})}))
 \to M_0(r,n)$ is finite and all fibers of
$\bM(r,k,\widehat{n}) \to Y$ are connected.
Then by (*), $Y \to M_0(r,n)$ is birational (at least over a field of
characteristic 0), the normality of $M_0(r,n)$ implies 
$Y \cong  M_0(r,n)$.
\end{NB}

We next prove (2).
Assume that $k=0$. Then the Poincar\'{e} dual of $\mu(C)$
is the closed subset of $\bM(r,0,\widehat{n})$ (see \cite{Bryan:1997}):
\begin{equation*}
\{ (E,\Phi) \in \bM(r,0,\widehat{n})|
E_{|C} \not \cong {\cal O}_C^{\oplus r} \}.
\end{equation*}
Hence we have a $\hT$-equivariant homomorphism 
${\cal O}_{\bM(r,0,\widehat{n})} \to {\cal O}_{\bM(r,0,\widehat{n})}(\mu(C))$,
which induces 
$\hT$-equivariant inclusions
${\cal O}_{\bM(r,0,\widehat{n})} \subset 
{\cal O}_{\bM(r,0,\widehat{n})}(d\mu(C))$
%\end{equation}
for $d \geq 1$.
By taking the direct images, we have 
a $\hT$-equivariant inclusion 
\begin{equation}\label{eq:direct}
{\cal O}_{M_0(r,n)}= 
\widehat{\pi}_*({\cal O}_{\bM(r,0,\widehat{n})})
 \subset \widehat{\pi}_*({\cal O}_{\bM(r,0,\widehat{n})}(d\mu(C))).
\end{equation}
We note that ${\cal O}_{M_0^{\operatorname{reg}}(r,n)} \to
\widehat{\pi}_*(
{\cal O}_{\bM(r,0,\widehat{n})}(d\mu(C)))_{|M_0^{\operatorname{reg}}(r,n)}$ is an 
isomorphism.
Since $M_0(r,n)$ is normal and
$\dim(M_0(r,n) \setminus M_0^{\operatorname{reg}}(r,n)) \leq \dim M_0(r,n)-2$,
the torsion freeness of 
$\widehat{\pi}_*({\cal O}_{\bM(r,0,\widehat{n})}(d\mu(C)))$
implies that
$\widehat{\pi}_*({\cal O}_{\bM(r,0,\widehat{n})}(d\mu(C)))={\cal O}_{M_0(r,n)}$
as a coherent sheaf.
Then \eqref{eq:direct} implies that 
\begin{equation*}
\widehat{\pi}_*({\cal O}_{\bM(r,0,\widehat{n})}(d\mu(C))) \cong
\widehat{\pi}_*({\cal O}_{\bM(r,0,\widehat{n})}) \cong {\cal O}_{M_0(r,n)}
\end{equation*} 
as $\hT$-equivariant sheaves. 
Thus (2) holds for $k=0$.
We next assume that $0<k<r$.
As we will see in Proposition \ref{nef&big},
$-\mu(C)$ is $\widehat{\pi}$-big.
Since $\widehat{\pi}$ is a $Gr(r,k)$-bundle over
$M_0^{\mathrm{reg}}(r,n)$ with
$\dim Gr(r,k)=k(r-k)>0$, 
$\widehat{\pi}_*({\cal O}_{\bM(r,0,\widehat{n})}(d\mu(C)))_{|M_0^{\mathrm{reg}}(r,n)}
=0$. Since $\widehat{\pi}_*({\cal O}_{\bM(r,0,\widehat{n})}(d\mu(C)))$ 
is torsion free,
(2) holds also in these cases.
\qed

Proof of Proposition~\ref{prop:vanish}:
We note that
we do not need the $\hT$-structure to prove the claim.
So we forget the $\hT$-action.
In order to apply the Kawamata-Viehweg vanishing theorem to the $\Q$-line
bundle $\mu(C)$, we first compute the canonical line bundle of
$\bM(r,k,\widehat{n})$. For a later use, 
we compute it as a $\hT$-equivariant sheaf.
\begin{Lemma}\label{lem:canonical}
We have the following equalities:
\begin{equation*}
\begin{split}
K_{M(r,n)} &= 
{\cal O}_{M(r,n)}(r\mu(K_{{\Bbb P}^2}+2\linf)) =
(t_1 t_2)^{-rn}{\cal O}_{M(r,n)} \in \Pic^{\hT}(M(r,n)) \otimes {\Bbb Q},\\
K_{\bM(r,k,\widehat{n})} & =
{\cal O}_{\bM(r,k,\widehat{n})}(r\mu(K_{\widehat{\Bbb P}^2}+2\linf))
= (t_1 t_2)^{-r\widehat{n}}
{\cal O}_{\bM(r,k,\widehat{n})}(r\mu(C)) \in
\Pic^{\hT}(\bM(r,k,\widehat{n})) \otimes {\Bbb Q}.
\end{split}
\end{equation*}
\end{Lemma}

\begin{proof}
We only compute $K_{\bM(r,k,\widehat{n})}$.
The computation of $K_{M(r,n)}$ is similar and simpler.
Let $({\cal E},\Phi)$ be the universal family on
$\widehat{\Bbb P}^2 \times \bM(r,k,\widehat{n})$
and $p_2\colon\widehat{\Bbb P}^2 \times \bM(r,k,\widehat{n})\to
\bM(r,k,\widehat{n})$ be the projection.
Note that $T_{\bM(r,k,\widehat{n})} \cong 
\Ext^1_{p_2}({\cal E},{\cal E}(-\linf))$.
Using the (equivariant) Grothendieck-Riemann-Roch theorem,
we see that 
{\allowdisplaybreaks
\begin{equation*}
\begin{split}
&c_1(K_{\bM(r,k,\widehat{n})})\\
 =\; & 
c_1\left( {\bf R}p_{2*}({\cal E}^{\vee} \otimes
{\cal E}(-\linf)) \right)\\
=\; &\left[p_{2*}\left( \ch({\cal E}^{\vee} \otimes {\cal E})
e^{-\linf} \Todd_{\widehat{\Bbb P}^2}\right)\right]_1\\
=\; &\left[ p_{2*}\left((r^2-r\Delta({\cal E})+
\ch_4({\cal E}^{\vee} \otimes {\cal E})+\cdots)
\left(1-\frac{2\linf+K_{\widehat{\Bbb P}^2}}{2}+\cdots \right)
\right)\right]_1\\
=\; & r\Delta({\cal E})/(K_{\widehat{\Bbb P}^2}+2\linf)
=r\mu(K_{\widehat{\Bbb P}^2}+2\linf),
\end{split}
\end{equation*}
where} $[...]_1$ means the codimension 1 component in
the equivariant Chow group
$A^{\hT}_{2\widehat{n}r-1}(\bM(r,k,\widehat{n}))_{\Bbb Q}$.
Since $\Pic^{\hT}(\bM(r,k,\widehat{n})) \cong
A^{\hT}_{2\widehat{n}r-1}(\bM(r,k,\widehat{n}))$
(see \cite[Lemma 1.3]{Y} or \cite[Theorem 1]{EG}),
$K_{\bM(r,k,\widehat{n})}= 
{\cal O}_{\bM(r,k,\widehat{n})}(r\mu(K_{\widehat{\Bbb P}^2}+2\linf))
\in \Pic^{\hT}(\bM(r,k,\widehat{n})) \otimes {\Bbb Q}$.
Since ${\cal E}_{|\linf \times \bM(r,k,\widehat{n})}
\cong \bigoplus_{\alpha=1}^r
{\cal O}_{\linf \times \bM(r,k,\widehat{n})}e_\alpha$
as $\hT$-sheaves,
we have 
$p_{2*}(\Delta({\cal E}_{|\linf \times \bM(r,k,\widehat{n})}))=0$
in $A^1_{\hT}(\bM(r,k,\widehat{n}))_\Q$. Therefore
${\cal O}_{\bM(r,k,n)}(\mu(\linf))
\cong {\cal O}_{\bM(r,k,n)}$ as 
a $\hT$-line bundle. 
Since $K_{{\Bbb P}^2} \cong (t_1 t_2)^{-1}{\cal O}_{{\Bbb P}^2}(-3\linf)$
\begin{NB}
$K_{{\Bbb P}^2} \cong {\cal O}(-\linf-\ell_x-\ell_y)$.
\end{NB}
and $K_{\widehat{\Bbb P}^2} \cong p^*(K_{{\Bbb P}^2})(C)$,
we get
\begin{equation*}
 K_{\bM(r,k,\widehat{n})}= 
 (t_1 t_2)^{-r\widehat{n}}{\cal O}_{\bM(r,k,\widehat{n})}(r\mu(C))
\in \Pic^{\hT}(\bM(r,k,\widehat{n})) \otimes {\Bbb Q},
\end{equation*}
where we used $\Delta({\cal E})/[\widehat{\Bbb P}^2]=\widehat{n}$.
\end{proof}
By this lemma,
${\cal O}_{\bM(r,k,\widehat{n})}(d\mu(C)) =
K_{\bM(r,k,\widehat{n})}((d-r)\mu(C))
\in \Pic^{\hT}(\bM(r,k,\widehat{n})) \otimes {\Bbb Q}$.
By the following proposition,
$-\mu(C)$ on $\bM(r,k,\widehat{n})$
is $\widehat{\pi}$-nef and $\widehat{\pi}$-big.
Then the relative version of the Kawamata-Viehweg vanishing theorem
\cite[Theorem 1-2-3]{KMM}
implies that 
\begin{equation*}
R^i \widehat{\pi}_*({\cal O}_{\bM(r,k,\widehat{n})}(d\mu(C))) 
\cong 
R^i \widehat{\pi}_*(K_{\bM(r,k,\widehat{n})}((d-r)\mu(C)+D))=0
\end{equation*}
for $d<r$, where $D=0$ in 
$A_{2\widehat{n}r-1}(\bM(r,k,\widehat{n})) \otimes {\Bbb Q}$.
Thus Proposition \ref{prop:vanish} holds.
\qed

\begin{Proposition}\label{nef&big}
The ${\Bbb Q}$-divisor $-\mu(C)$ on $\bM(r,k,\widehat{n})$
is $\widehat{\pi}$-nef and $\widehat{\pi}$-big.
\end{Proposition}
Before proving this proposition, we shall treat the remaining case
(3) $d=r$ and 
finish the proof of
Proposition \ref{prop:key}.
%\begin{Remark}
%If $k=0$, then $\widehat{\pi}$ is birational, which implies that
%all ${\Bbb Q}$-divisors are $\widehat{\pi}$-big.
%\end{Remark}
Since 
$\widehat{\pi}_{*}({\cal O}_{\bM(r,k,\widehat{n})})
={\cal O}_{M_0(r,n)}$ in $K^{\hT}(M_0(r,n))$,
the Grothendieck-Serre duality implies that
$\widehat{\pi}_{*}(K_{\bM(r,k,\widehat{n})})
=(-1)^{k(r-k)}K_{M_0(r,n)}=(-1)^{k(r-k)}\pi_*(K_{M(r,n)})$
in $K^{\hT}(M_0(r,n))$.
\begin{NB}
$M_0(r,n)$ has rational singularities.
Hence it is Cohen-Macauley. In particular the dualizing complex 
is represented by ${\cal O}_{M_0(r,n)}$.
\end{NB}
By using Lemma \ref{lem:canonical}, we get (3).

%\subsection{Proof of the $\widehat{\pi}$-nefness and $\widehat{\pi}$-bigness 
%of $-\mu(C)$.}

\begin{NB}
I have re-written the following paragraph. Please check it. (May 20) 
\end{NB}

The proof of \propref{nef&big} occupies the rest of this section.
Let us briefly sketch the idea of the proof.
Let $M(r,-H,n')$ be the moduli space of $H$-stable sheaves $E$ on
$\proj^2$ with $c_1(E) = c_1$, $\Delta(E) = n'= n+(r-1)(2r-1)/r$. The
first step of the proof is to construct an embedding $M(r,n)\to
M(r,-H,n')$ by using idea in \cite{parabolic}. Next we show that there
is an induced commutative diagram
\begin{equation}\label{eq:diagram1}
\begin{matrix}
M(r,n) & \hookrightarrow & M(r,-H,n')\\
{\scriptstyle {\pi}}\downarrow {\;}& & 
{\;}\downarrow \scriptstyle{{\pi}'}\\
M_0(r,n) & \to & M_0(r,-H,n'),
\end{matrix}
\end{equation}
where $M_0(r,-H,n')$ is the Uhlenbeck compactification of the open
subset consisting of locally free sheaves in $M(r,-H,n')$, and $\pi'$
is the morphism defined exactly as $\pi$.
Similarly we can make a commutative diagram on blowup:
\begin{equation}\label{eq:cd}
\begin{matrix}
\bM(r,k,\widehat{n}) & \hookrightarrow & \bM(r,kC-H,\widehat{n}')\\
{\scriptstyle \widehat{\pi}}\downarrow {\;}& & 
{\;}\downarrow \scriptstyle{\widehat{\pi}'}\\
M_0(r,n) & \to & M_0(r,-H,n'),
\end{matrix}
\end{equation}
where $\bM(r,kC-H,\widehat{n}')$ is the moduli space of $(H-\ve
C)$-stable sheaves on $\bp$ with $\widehat{n}'-\widehat{n}=n'-n$.
 From the diagram it is enough to show that $-\mu(C)$ is
$\widehat{\pi}'$-nef and $\widehat{\pi}'$-big. But this assertion
follows from the following known results:
1) The divisor $\mu(H-\ve C)$ is nef and big as it gives a
birational morphism $\bM(r,kC-H,\widehat{n}') \to
\bM_0(r,kC-H,\widehat{n}')$ \cite{Li:2}.
\begin{NB}
Please check that the target of the birational morphism is correct.
(May 20)
\end{NB}
2) $\mu(H)$ is the pull-back of a ${\Bbb Q}$-Cartier divisor ${\cal
  H}$ on $M_0(r,-H,n')$ by the construction of $\widehat{\pi}'$ in
\cite[Appendix~F]{lecture}.

\subsection{An embedding $M(r,n) \to M(r,-H,n')$}
\begin{NB}
The embedding is not $\hT$-equivariant.
\end{NB}
We construct an embedding $M(r,n) \to M(r,-H,n')$ using parabolic
sheaves \cite{M-Y}.
Our construction is inspired by \cite{parabolic}.
We collect basic facts on parabolic sheaves in \cite{M-Y}. 
\begin{Definition}
A {\it parabolic sheaf} $(E_*,\alpha_*)$ with respect to $\linf$ is 
a pair of a filtration 
\begin{equation*}
E_*:E(-\linf) \subset E_l \subset \dots \subset E_2 
\subset E_1=E
\end{equation*}
on ${\Bbb P}^2$
and a sequence of rational numbers
$0 \leq \alpha_1 <\alpha_2<\dots<\alpha_l<1$.

\end{Definition}
In this paper, we consider parabolic sheaves 
$(E_*,\alpha_1,\alpha_2)$ such that
\begin{equation}
E_*:E(-\linf) \subset E' \subset E
\end{equation}
and
$E/E' = {\cal O}_{\linf}(r-1)$ in the Grothendieck group
$K(\linf)$.
\begin{Definition}
A parabolic sheaf 
$(E_*,\alpha_1,\alpha_2)$ is {\it $\mu$-semi-stable}, if
$E$ is torsion free and
\begin{equation}\label{eq:stability2}
\begin{split}
 & \frac{\deg G +\alpha_1 \rank_{\linf} (G/G \cap E')+
 \alpha_2 \rank_{\linf} (G \cap E'/G(-\linf))}{\rank G}\\
  \leq\; &
 \frac{\deg E +\alpha_1 \rank_{\linf} (E/E')+
 \alpha_2 \rank_{\linf} (E'/E(-\linf))}{\rank E}
 \end{split}
 \end{equation}
 for a saturated subsheaf $G \subset E$ with $0<\rank G<\rank E$,
 where $\deg$ denotes the degree with respect to $H$
  and $\rank_{\linf}$ the rank on $\linf$.
If the inequality is strict for all $G$, then
we say $(E_*,\alpha_1,\alpha_2)$ is {\it $\mu$-stable}.
\end{Definition}
We set $\alpha=\alpha_2-\alpha_1$, then $0<\alpha<1$ and
\eqref{eq:stability2} is equivalent to the following condition:
\begin{equation}\label{eq:stability}
 \frac{\deg G -\alpha \rank_{\linf} (G/G \cap E')}{\rank G} \leq
 \frac{\deg E -\alpha \rank_{\linf} (E/E')}{\rank E}
 \end{equation}
 for a saturated subsheaf $G \subset E$ with $0<\rank G<\rank E$.
Therefore we define our parabolic sheaf as a filtration
$E_*:E(-\linf) \subset E' \subset E$ with a parameter $\alpha$, and we
define the stability of $(E_*,\alpha)$ as the condition 
\eqref{eq:stability}.

 From now on, we fix the parameter $\alpha$ with
$0<\alpha <r/(r-1)^2$. Thus our parabolic sheaf is 
the filtration $E_*:E(-\linf) \subset E' \subset E$. 
\begin{Definition}
${\frak M}$ denotes the moduli space of $\mu$-stable parabolic sheaves
$(E_*, \alpha)$ on ${\Bbb P}^2$ such that
$(\rank E,c_1(E),\Delta(E))=(r,0,n)$
and $E/E'= {\cal O}_{\linf}(r-1)$ in 
$K(\linf)$.
\end{Definition}
An easy computation shows that
$(c_1(E'),\Delta(E'))=(-H,n')$.

\begin{Lemma}\label{lem:stability}
Let $E_*:E(-\linf) \subset E' \subset E$ be a parabolic sheaf
with $(\rank E,c_1(E),\Delta(E))=(r,0,n)$
and $E/E'= {\cal O}_{\linf}(r-1)$.
\begin{enumerate}
\item
$E_*$ is $\mu$-stable if and only if
$E$ is $\mu$-semi-stable and
$\rank_{\linf}(G/G \cap E')=1$ 
for all saturated subsheaf $0 \ne G \subset E$ with $\deg(G)=0$.
\item
$E_*$ is $\mu$-stable if and only if
$E'$ is $\mu$-stable.
Hence we have a morphism
$\phi:{\frak M} \to M(r,-H,n')$.
\end{enumerate}
\end{Lemma}

\begin{proof}

We first prove (1).
Assume that $E_*$ is $\mu$-stable.
If $E$ is not $\mu$-semi-stable, then 
there is a saturated subsheaf $G \subset E$ with $\deg G >0$.
Then
\begin{equation*}
\begin{split}
& \frac{\deg E -\alpha \rank_{\linf} (E/E')}{\rank E}
- \frac{\deg G -\alpha \rank_{\linf} (G/G \cap E')}{\rank G}\\ 
 \leq\; & \frac{-1}{r-1}+\alpha
\left(\frac{\rank_{\linf}(G/E' \cap G)}{\rank G}-\frac{1}{r}\right)\\
\leq\; & \frac{-1}{r-1}+\alpha \left(1-\frac{1}{r}\right) 
  \leq  \frac{1}{r(r-1)}((r-1)^2\alpha-r)<0,
  \end{split}
 \end{equation*}
which is a contradiction.
Therefore $E$ is $\mu$-semi-stable.
If $E$ is properly $\mu$-semi-stable, then
there is a saturated subsheaf $G \subset E$ with $\deg G=0$.
In this case, we must have
${\rank_{\linf}(G/G \cap E')}/{\rank G}>{1}/{r}$.
Since $G/G \cap E' \subset E/E'$ and $\rank_{\linf}(E/E')=1$,
$\rank_{\linf}(G/G \cap E')=1$.
We show the inverse direction.
We assume that $E$ is $\mu$-stable.
Then for a saturated subsheaf $G \subset E$, we have $\deg G < 0$.
Then we see that
\begin{equation*}
\begin{split}
& \frac{\deg E -\alpha \rank_{\linf} (E/E')}{\rank E}
- \frac{\deg G -\alpha \rank_{\linf} (G/G \cap E')}{\rank G}\\ 
 \geq\; & \frac{1}{r-1}
+\alpha
\left(\frac{\rank_{\linf}(G/E' \cap G)}{\rank G}-\frac{1}{r}\right)
 \geq \frac{1}{r-1}-\frac{\alpha}{r} >0.
%\geq -\frac{1}{r-1}((r-1)\alpha-1)>0.
 \end{split}
 \end{equation*}
Thus $E_*$ is $\mu$-stable.
If $E$ is properly $\mu$-semi-stable, we assume that  
$\rank_{\linf}(G/G \cap E')=1$ for all saturated subsheaf 
$G \subset E$ with $\deg G=0$.
Then we also see that $E_*$ is $\mu$-stable.

We next prove (2).
If $E'$ is not $\mu$-stable, then there is a saturated subsheaf
$G' \subset E'$ with $0<\rank G' <\rank  E'$ and 
\begin{equation*}
\frac{\deg G'}{\rank G'} \geq \frac{\deg(E')}{\rank E'}=-\frac{1}{r}.
\end{equation*}
Therefore $\deg G' \geq 0$.
Since $E'/G'$ is torsion free, $G'_{|\linf} \to E'_{|\linf}$ is injective,
and hence $G' \cap E'(-\linf)=G'(-\linf)$. 
We set $G:=(G' \cap E(-\linf))(\linf)$.
Since $(G' \cap E(-\linf)) \cap E'(-\linf)=
G' \cap E'(-\linf)=G'(-\linf)$,
we get $G \cap E'=G'$.
By (1), $E$ is $\mu$-semi-stable, therefore $\deg G \leq 0$.
As $0 \leq \deg G' \leq \deg G \leq 0$, we get $\deg G=\deg G'=0$.
Hence $\rank_{\linf} (G/G')=0$, which implies that
$E(-\ell_\infty) \subset E' \subset E$ is not $\mu$-semi-stable. 

Conversely suppose $E'$ is $\mu$-stable.
Then for a subsheaf $G \subset E$,
we have $\deg G \cap E'<0$.
As $\deg G=\deg(G \cap E')+\rank_{\linf}(G/G \cap E')$,
we have $\deg G < 0$ or $\deg G=0$ and $\rank(G/G \cap E')=1$.
By (1), $E_*$ is $\mu$-stable.
\end{proof}

Let ${\frak M}'$ be the open subscheme of ${\frak M}$ such that
 \begin{equation*}
 \begin{split}
 {\frak M}'=&
 \left\{E_* \in {\frak M} \,\left|
  \text{ $E/E'$ and $E'/E(-\linf)$ are }
  \text{ semi-stable locally free sheaves on $\linf$ } 
 \right. \right\}\\
 = &
 \{E_* \in {\frak M}\,|\text{ $E/E' \cong {\cal O}_{\linf}(r-1)$},
 \text{ $E'/E(-\linf) \cong {\cal O}_{\linf}(-1)^{\oplus (r-1)}$} \}.
 \end{split}
 \end{equation*}
 If $E(-\linf) \subset E' \subset E$ belongs to ${\frak M}'$,
 then $E_{|\linf}$ is a locally free ${\cal O}_{\linf}$-module.
 Hence $E$ is locally free in a neighborhood of $\linf$, and hence
 $E'=\ker(E \to {\cal O}_{\linf}(r-1))$ is also locally free 
 in a neighborhood of $\linf$.
 
\begin{Lemma}\label{lem:M'}
Let $E_*:E(-\linf) \subset E' \subset E$ be a point of ${\frak M}$.
\begin{enumerate}
\item
If $E_*$ belongs to ${\frak M}'$, then
$E'_{|\linf} \cong 
{\cal O}_{\linf}(r-2) \oplus {\cal O}_{\linf}(-1)^{\oplus (r-1)}$.
\item
If $E'_{|\ell_{\infty}} \cong 
{\cal O}_{\linf}(r-2) \oplus {\cal O}_{\linf}(-1)^{\oplus (r-1)}$, 
then 
\begin{equation}\label{eq:E}
E=\ker(E' \to \Hom(E',{\cal O}_{\linf}(-1))^{\vee} 
\otimes {\cal O}_{\linf}(-1)) 
\otimes {\cal O}_{{\Bbb P}^2}(\linf).
\end{equation}
In particular $E_* \in {\frak M}$ is uniquely determined by $E'$.
Moreover $E_*$ belongs to ${\frak M}'$.
%$E$ is locally free along $\linf$,
%$E/E' \cong {\cal O}_{\linf}(r-1)$ and 
%$E_{|\linf}$ fits in an exact sequence
%\begin{equation}
%0 \to {\cal O}_{\linf}(-1)^{\oplus (r-1)} \to E_{|\linf}
%\to {\cal O}_{\linf}(r-1) \to 0.
%\end{equation}
\end{enumerate}
\end{Lemma}

\begin{proof}
(1)
By the filtration
$E'(-\linf) \subset E(-\linf) \subset E' \subset E$ induced by $E_*$,
we have an exact sequence
\begin{equation*}
0 \to (E/E') \otimes {\cal O}_{{\Bbb P}^2}(-\linf) \to
E'_{|\linf} \to E'/E(-\linf) \to 0.
\end{equation*}
Then this exact sequence splits and we get our claim.

(2)
We set $L':=E'/E(-\linf)$. Then
$L'$ is an ${\cal O}_{\linf}$-module of rank $r-1$ with 
$\deg L'=-(r-1)$.
Let $T$ be the torsion submodule of $L'$ and
\begin{equation*}
 0 \subset F_1 \subset F_2 \subset \dots \subset F_s=L'/T
\end{equation*}
the Harder-Narasimhan filtration of $L'/T$.
We set $F_i/F_{i-1} \cong {\cal O}_{\linf}(a_i)^{\oplus n_i}$.
Then 
$a_1>a_2>\cdots>a_s$ with 
$\sum_{i=1}^s n_i a_i+\deg(T)=\deg L'=-(r-1)$.
Since we have a surjective homomorphism 
\begin{equation*}
 E'_{|\linf} \to L' \to {\cal O}_{\linf}(a_s),
\end{equation*}
we get 
$a_s \geq -1$.
Then we see that $T=0$, $s=1$ and
$L' \cong {\cal O}_{\linf}(-1)^{\oplus (r-1)}$. 
Moreover we have a commutative diagram
\begin{equation*}
\begin{CD}
E' @>>> L'\\
@| @VV{\xi}V \\
E' @>>> \Hom(E,{\cal O}_{\linf}(-1))^{\vee} \otimes {\cal O}_{\linf}(-1),
\end{CD}
\end{equation*}
where $\xi$ is an isomorphism.
Hence \eqref{eq:E} holds.
By the filtration
$E'(-\linf) \subset E(-\linf) \subset E'$ induced by $E_*$,
we see that $E/E' \cong {\cal O}_{\linf}(r-1)$. Therefore
(2) holds.
\end{proof}

We set 
\begin{equation*}
S:=\{E' \in M(r,-H,n')|E'_{|\linf} \cong 
{\cal O}_{\linf}(r-2) \oplus {\cal O}_{\linf}(-1)^{\oplus (r-1)} \}.
\end{equation*}
$S$ is a locally closed subscheme of $M(r,-H,n')$, 
where we use the reduced scheme structure on $S$.
By Lemma \ref{lem:M'} (1), we get $\phi({\frak M}') \subset S$.
\begin{Lemma}\label{lem:isom1}
${\frak M}' \to S$ is an isomorphism.
\end{Lemma}

\begin{proof}
Let ${\cal E}'$ be the universal family on
$S \times {\Bbb P}^2$ and
$q_S:S \times {\Bbb P}^2 \to S$ the projection.
By the definition of $S$,
$\Hom({\cal E}'_{|\{s \} \times \linf},{\cal O}_{\linf}(-1))
\cong {\Bbb C}^{\oplus (r-1)}$ for $s \in S$.
Since $S$ is reduced, the base change theorem implies that
$U:=\Hom_{q_S}({\cal E}'_{|S \times \linf},
{\cal O}_{S \times \linf}(-1))$ is a locally free sheaf on $S$ and
we have a family of homomorphisms
\begin{equation*}
f:{\cal E}'_{|S \times {\Bbb P}^2} \to U^{\vee} \boxtimes {\cal O}_{\linf}(-1).
\end{equation*}
We set ${\cal E}:=(\ker f)(\linf)$.
Then we have a family of parabolic sheaves
${\cal E}(-\linf) \subset {\cal E}' \subset {\cal E}$.
They are $\mu$-stable by Lemma \ref{lem:stability} (2).
Thus we have a morphism $\psi:S \to {\frak M}'$.
It is easy to see that $\phi \circ \psi$ is the identity.
By Lemma \ref{lem:M'} (2), we also see that
$\psi \circ \phi$ is the identity.
Therefore $\phi$ is an isomorphism. 
\end{proof}
Let $({\cal E},\Phi)$ be the universal family on 
$M(r,n) \times {\Bbb P}^2$.
Then we have a surjective homomorphism
\begin{equation*}
\Psi:{\cal E} \to {\cal O}_{M(r,n) 
\times \linf}^{\oplus r} \overset{\phi}{\to}
{\cal O}_{M(r,n) \times \linf}(r-1)
\end{equation*}
where $\phi:=(z_1^{r-1},z_1^{r-2}z_2,...,z_2^{r-1})$.
We set ${\cal E}':=\ker \Psi$. Then
we have a family of parabolic sheaves
${\cal E}(-\linf) \subset {\cal E}' \subset {\cal E}$.

\begin{Lemma}
${\cal E}(-\linf) \subset {\cal E}' \subset {\cal E}$ is a family of
$\mu$-stable parabolic sheaves.
\end{Lemma}

\begin{proof}
For $s \in M(r,n)$, we set
$E:={\cal E}_{|\{s \} \times {\Bbb P}^2}$ and
$E':={\cal E}'_{|\{s \} \times {\Bbb P}^2}$. 
Since $E_{\linf} \cong {\cal O}_{\linf}^{\oplus r}$, $E$ is 
$\mu$-semi-stable.
Assume that there is a saturated proper subsheaf $G$ of $E$ with $\deg G=0$.
Then, since $E/G$ is torsion free,
$G_{|\linf}$ is a subsheaf of $E_{|\linf} \cong {\cal O}_{\linf}^{\oplus r}$.
Hence $G_{|\linf}$ is a direct summand of $E_{|\linf}$, which implies that
$\rank_{\linf}(G/G \cap E')=1$.
By Lemma \ref{lem:stability} (1),
$E(-\linf) \subset E' \subset E$ is $\mu$-stable.
\end{proof}
Hence we have a morphism
$M(r,n) \to {\frak M}$.
Let ${\frak M}''$ be the open subscheme of 
${\frak M}'$ such that $E(-\linf) \subset E' \subset E$
satisfies $E_{|\linf} \cong {\cal O}_{\linf}^{\oplus r}$. 

\begin{Lemma}\label{lem:isom2}
We have an isomorphism $M(r,n) \to {\frak M}''$. 
\end{Lemma}

\begin{proof}
We construct the inverse map.
For the family of parabolic sheaves
${\cal E}(-\linf) \subset {\cal E}' \subset {\cal E}$
on ${\frak M}'' \times {\Bbb P}^2$, we see that
${\cal E}/{\cal E}' \cong L \boxtimes {\cal O}_{\linf}(r-1)$
for a line bundle $L$ on ${\frak M}''$.
By taking the tensor product with 
$q_{{\frak M}''}^*(L^{\vee})$ to the family, 
we may assume that 
${\cal E}/{\cal E}' \cong 
{\cal O}_{{\frak M}''} \boxtimes {\cal O}_{\linf}(r-1)$.
We set
$g:{\cal E} \to 
{\cal O}_{{\frak M}''} \boxtimes {\cal O}_{\linf}(r-1)$.
Since ${\cal E}_{|\linf}$ is a family of ${\cal O}_{\linf}^{\oplus r}$, 
$g$ induces an isomorphism
$q_{{\frak M}''*}({\cal E}_{|\linf}) \to 
{\cal O}_{{\frak M}''} \otimes H^0(\linf,{\cal O}_{\linf}(r-1))$.
Then we have a commutative diagram
\begin{equation*}
\begin{CD}
 q_{{\frak M}''*}({\cal E}_{|\linf}) \boxtimes {\cal O}_{\linf}
@>{\sim}>> {\cal O}_{{\frak M}''} \boxtimes {\cal O}_{\linf}^{\oplus r}\\
@VVV @VV{\phi}V\\
{\cal E}_{|\linf} @>{g_{|\linf}}>> {\cal O}_{{\frak M}'' \times \linf}(r-1)
\end{CD}
\end{equation*} 
Since the left vertical map is an isomorphism,
we can find a homomorphism 
$\Phi:{\cal E} \to {\cal O}_{{\frak M}'' \times \linf}^{\oplus r}$
such that $\phi \circ \Phi=g$. 
\end{proof}
By Lemma \ref{lem:isom1} and  Lemma \ref{lem:isom2},
we have a sequence of morphisms:
\begin{equation*}
M(r,n) \overset{\sim}{\to} {\frak M}'' \hookrightarrow
{\frak M}' \overset{\sim}{\to} S \hookrightarrow M(r,-H,n').
\end{equation*}
Thus we obtain the following proposition.
\begin{Proposition}
We have an immersion
$\iota:M(r,n) \to M(r,-H,n')$.
\end{Proposition}

\subsection{Construction of $M_0(r,n) \to M_0(r,-H,n')$}

Let $\pi':M(r,-H,n') \to M_0(r,-H,n')$ be the contraction map to the
Uhlenbeck's compactification of the open subset consisting of locally
free sheaves in $M(r,-H,n')$. Since $\pi'$ is surjective, replacing
$M_0(r,-H,n')$ by its normalization, we may assume that $M_0(r,-H,n')$
is normal.  Then by Corollary~\ref{cor:equivalence}, $M_0(r,-H,n')$ is
the quotient of the equivalence relation which defines the contraction
to the Uhlenbeck compactification.  We consider the composition
$\varpi:M(r,n) \to M(r,-H,n') \to M_0(r,-H,n')$.  Then

\begin{Lemma}
 $\varpi((E_1,\Phi_1))=\varpi((E_2,\Phi_2))$
if and only if $(E_1^{\vee \vee},\Phi_1) \cong (E_2^{\vee \vee},\Phi_2)$
and $\Supp(E_1^{\vee \vee}/E_1)=\Supp(E_2^{\vee \vee}/E_2)$.
\end{Lemma}

\begin{proof}
We set $E_i':=\iota((E_i,\Phi_i))$, $i=1,2$. 
Then $\varpi((E_1,\Phi_1))=\varpi((E_2,\Phi_2))$ if and only if
$(E_1')^{\vee \vee} \cong (E_2')^{\vee \vee}$
and $\Supp((E_1')^{\vee \vee}/E_1')=\Supp((E_2')^{\vee \vee}/E_2')$.
Since $E_i'$, $i=1,2$ fits in an exact sequence
\begin{equation*}
0 \to E_i' \to E_i \to {\cal O}_{\linf}(r-1) \to 0
\end{equation*}
and $E_i$ is locally free along $\linf$,
we have an exact sequence
\begin{equation*}
0 \to (E_i')^{\vee \vee} \to E_i^{\vee \vee} \to {\cal O}_{\linf}(r-1) \to 0.
\end{equation*}
Moreover $E_i^{\vee \vee} \to {\cal O}_{\linf}(r-1)$ is uniquely determined by
$(E_i')^{\vee \vee}$.
Therefore our lemma holds.
\end{proof}
We thus get the diagram \eqref{eq:diagram1} except the bottom
arrow. Now the morphism $M_0(r,n) \to M_0(r,-H,n')$ exists as
$M_0(r,n)$ is normal by \lemref{lem:morphism}.
\begin{NB}
Please check my change \corref{cor:equivalence} $\to$ \lemref{lem:morphism}
is correct. (May 20) 
\end{NB}

\subsection{Proof of Proposition \ref{nef&big}}

In the same way, we have an immersion
$\bM(r,k,\widehat{n}) \to \bM(r,kC-H,\widehat{n}')$ and the 
commutative diagram \eqref{eq:cd}.
% \begin{equation}
% \begin{matrix}                                                                  
% \bM(r,k,\widehat{n}) & \hookrightarrow & \bM(r,kC-H,\widehat{n}')\\
% {\scriptstyle \widehat{\pi}}\downarrow {\;}& & 
% {\;}\downarrow \scriptstyle{\widehat{\pi}'}\\
% M_0(r,n) & \to & M_0(r,-H,n')
% \end{matrix}
% \end{equation}
%
%where $\widehat{n}'-\widehat{n}=n'-n$.

As we already mentioned before, this implies that $-\mu(C)$ is
$\widehat{\pi}'$-nef and $\widehat{\pi}'$-big. Here $\mu$ is defined
as $\mu(\bullet):=\Delta({\cal E}')/[\bullet]$ by the universal family
${\cal E}'$ on $\widehat{\Bbb P}^2  \times \bM(r,kC-H,\widehat{n}')$.
This is enough for our purpose as two universal sheaves $\cal E$ and
$\cal E'$ are isomorphic outside $\linf$, so $\mu(C)$ turns out to be
the same for $\cal E$ and $\cal E'$.
Thus we complete the proof of Proposition \ref{nef&big}.
\qed

%\section{K-theoretic partition function}

\section{Behavior at $\ve_1,\ve_2 = 0$}

We prove a part of Nekrasov's conjecture and derive an equation
satisfied by the `genus $0$' part of the partition function in this
section.

\begin{NB}
I use the inverse of the pull-back to define the action
of $\hT$ on coherent sheaves on $M(r,n)$.
So in order to match the notations in \cite{part1}, we need to
change $t_i$ by $t_i^{-1}$ and $e_\alpha$
by $e_\alpha^{-1}$ in \eqref{eq:l} and \eqref{eq:blow-up1}.

\subsection{$\hT$-actions on the canonical line bundles}
Let $(E,\Phi)$ be a fixed point of $\hT$-action corresponding to
$\vec{Y} = (Y_1,\dots, Y_r)$. 
Then the $\hT$-module structure of
$\det T_{(E,\Phi)} M(r,n)$ is given by
\begin{equation*}
   \det T_{(E,\Phi)} M(r,n) 
   =(t_1 t_2)^{rc_2(E)}.
\end{equation*}

Let $(E,\Phi)$ be a fixed point of $\hT$-action corresponding to
$(\vec{k}, \vec{Y}^1,\vec{Y}^2)$. Then the $\hT$-module structure of
$\det T_{(E,\Phi)} \bM(r,k,\widehat{n})$ is given by
\begin{equation*}
\begin{split}
   \det T_{(E,\Phi)} \bM(r,k,\widehat{n})
   &=t_2^{rc_2(\vec{Y}^1)}t_1^{rc_2(\vec{Y}^2)}
 \prod_{\alpha < \beta}
\left(\frac{e_\beta}{e_\alpha}\right)^{k_\alpha-k_\beta}\\
&=\left(t_1 t_2 \right)^{r\Delta(E)}
t_1^{-rc_2(\vec{Y}^1)}t_2^{-rc_2(\vec{Y}^2)}
\prod_{\alpha < \beta}
\left(t_1 t_2 \right)^{-\frac{(k_\alpha-k_\beta)^2}{2}}
\left(\frac{e_\beta}{e_\alpha}\right)^{k_\alpha-k_\beta} \\
&=\left(t_1 t_2 \right)^{r(\Delta(E)-\frac{k(r-k)}{2r})}
\left(t_1 t_2 \right)^{\frac{k(r-k)}{2}}
{\cal O}_{\bM(r,k,\widehat{n})}(-r\mu(C))_{|(E,\Phi)}.
\end{split}
\end{equation*}
Thus 
\begin{equation*}
K_{\bM(r,k,\widehat{n})}=(t_1 t_2)^{r \widehat{n}}
{\cal O}_{\bM(r,k,\widehat{n})}(-r\mu(C))
\end{equation*}
in $K^{\hT}(\bM(r,k,\widehat{n})) \otimes_{R(\hT)} {\cal R}$.
\end{NB}

\subsection{Nekrasov's conjecture}

We collect some properties of $l^{\vec{k}}_{\vec{\alpha}}$.
\begin{Lemma}\label{lem:lsym}
\textup{(1)}
\(
   l^{\vec{k}}_{\vec{\alpha}}(\ve_1,\ve_2,\vec{a}) = 
   l^{\vec{k}}_{\vec{\alpha}}(\ve_2,\ve_1,\vec{a}).
\)

\textup{(2)}
\(
   l^{\vec{k}}_{\vec{\alpha}}(\ve_1,\ve_2,\vec{a}) = 
   (-e^{-\bbeta\langle\vec{a},\alpha\rangle})
   ^{\langle\vec{k},\alpha\rangle(\langle\vec{k},\alpha\rangle-1)/2}\,
   (e^{\bbeta(\ve_1+\ve_2)})^{{\langle\vec{k},\alpha\rangle
   (\langle\vec{k},\alpha\rangle^2-1)}/{6}}\,
   l^{-\vec{k}}_{-\vec{\alpha}}(-\ve_1,-\ve_2,\vec{a}).
\)

\textup{(3)}
$l^{\vec{k}}_{\vec{\alpha}}(\ve_1,\ve_2,\vec{a})$ is regular at
$(\ve_1,\ve_2) = 0$ and
\begin{equation*}
   l^{\vec{k}}_{\vec{\alpha}}(0,0,\vec{a}) = 
    (1-e^{-\bbeta\langle \vec{a},\vec\alpha\rangle})^{
     \langle\vec{k},\alpha\rangle(\langle\vec{k},\vec\alpha\rangle-1)/2}.
\end{equation*}
\end{Lemma}

Since
\begin{equation*}
\prod_{\vec\alpha \in \Delta}(-e^{-\bbeta\langle\vec{a},\vec\alpha\rangle})
   ^{\langle\vec{k},\vec\alpha\rangle(\langle\vec{k},\vec\alpha\rangle-1)/2}
=\prod_{\alpha \in \Delta_+}
e^{\bbeta\langle\vec{a},\vec\alpha\rangle \langle\vec{k},\vec\alpha\rangle}
(-1)^{\langle\vec{k},\vec\alpha\rangle}
=e^{r\bbeta(\vec{k},\vec{a})}(-1)^{2\langle \vec{k},\rho \rangle}
=e^{r\bbeta(\vec{k},\vec{a})},
\end{equation*}
we get
\begin{equation}\label{eq:l-sym}
\prod_{\vec\alpha \in \Delta}
l^{\vec{k}}_{\vec{\alpha}}(\ve_1,\ve_2,\vec{a})=
e^{r\bbeta(\vec{k},\vec{a})}\prod_{\vec\alpha \in \Delta}
l^{-\vec{k}}_{-\vec\alpha}(-\ve_1,-\ve_2,\vec{a}).
\end{equation}

The following will be used only afterwards, but it illustrates a usage
of the blowup equations in \thmref{thm:key2}.

\begin{Lemma}\label{lem:blow-up}
\begin{equation*}
\Zin(\ve_1,-2\ve_1,\vec{a};\q, \bbeta)
  =\Zin(2\ve_1,-\ve_1,\vec{a};\q,\bbeta). 
\end{equation*}
\end{Lemma}

\begin{proof}
We apply the relations \eqref{eq:blow-up1} and \eqref{eq:blow-up2} for
$d=r,0$ after setting $\ve_2 = -\ve_1$. We omit $\bbeta$ by letting it
be $1$ for brevity. Then
\begin{equation*}
\begin{split}
  \sum_{\vec{k}}
 &  \frac{\q^{(\vec{k},\vec{k})/2}}
   {\prod_{\vec\alpha\in\Delta} l^{\vec{k}}_{\vec{\alpha}}(\ve_1,-\ve_1,\vec{a})}
   \left(e^{(\vec{k},\vec{a})r}
     \Zin(\ve_1,-2\ve_1,\vec{a}+\ve_1\vec{k};t_1^{r/2}\q) 
     \Zin(2\ve_1,-\ve_1,\vec{a}-\ve_1\vec{k};t_1^{-r/2}\q) \right.\\
& \quad - \left.
   \Zin(\ve_1,-2\ve_1,\vec{a}+\ve_1\vec{k};t_1^{-r/2}\q)
   \Zin(2\ve_1,-\ve_1,\vec{a}-\ve_1\vec{k};t_1^{r/2}\q)\right)=0.
\end{split}
\end{equation*}
%We note that
%\begin{equation}
%\frac{e^{r\frac{(\vec{k},\vec{a})}{2}}}
%{\prod_{\vec\alpha\in\Delta} l^{\vec{k}}_{\vec{\alpha}}(\ve_1,-\ve_1,\vec{a})}
%=\frac{e^{\frac{r(-\vec{k},\vec{a})}{2}}}
%{\prod_{\vec\alpha\in\Delta} l^{-\vec{k}}_{\alpha}(\ve_1,-\ve_1,\vec{a})}
%\end{equation}
%and
%\begin{equation}
%\begin{split}
%&
% \Zin(\ve_1,-2\ve_1,\vec{a}+\ve_1\vec{k};t_1^d\q)
%   \Zin(2\ve_1,-\ve_1,\vec{a}-\ve_1\vec{k};t_1^{-d}\q)\\
%=&
%\sum_{n,m} \Zin_n(\ve_1,-2\ve_1,\vec{a}+\ve_1\vec{k})(t_1^{d-\frac{r}{2}}\q)^n
%   \Zin_m(2\ve_1,-\ve_1,\vec{a}-\ve_1\vec{k})(t_1^{-d+\frac{r}{2}}\q)^m.
%\end{split}
%\end{equation}

Let us expand $\Zin$ as in \eqref{eq:instantonpart}.
Then the above equation implies
\begin{equation*}
\begin{split}
 & \left(
   Z_n(\ve_1,-2\ve_1,\vec{a})
   -Z_n(2\ve_1,-\ve_1,\vec{a}) \right)
(t_1^{{rn}/{2}}-t_1^{-{rn}/{2}})\\
= &
-\sum_{\substack{(\vec{k},\vec{k})/2+l+m = n\\
      l\neq n, m\neq n}}
\begin{aligned}[t]
  & \frac{e^{r{(\vec{k},\vec{a})}/{2}}
   Z_m(\ve_1,-2\ve_1,\vec{a}+\ve_1\vec{k})
   Z_l(2\ve_1,-\ve_1,\vec{a}-\ve_1\vec{k})}
  {\prod_{\vec\alpha\in\Delta} l^{\vec{k}}_{\vec{\alpha}}(\ve_1,-\ve_1,\vec{a})}
  \\
  & \qquad\qquad
  \times\left(e^{r{(\vec{k},\vec{a})}/{2}}
    t_1^{{r(m-l)}/{2}}-
    e^{{r(-\vec{k},\vec{a})}/{2}}t_1^{{-r(m-l)}/{2}}\right). 
\end{aligned}
 \end{split}
\end{equation*}
Let us show that
\(
Z_n(\ve_1,-2\ve_1,\vec{a})
   = Z_n(2\ve_1,-\ve_1,\vec{a})
\)
by using the induction on $n$. It holds for $n = 0$ as $Z_0 = 1$.
Suppose that it is true for $l,m< n$. Then the right hand side of the
above equation vanishes, as terms with $(\vec{k},l,m)$ and $(-\vec{k},
m, l)$ cancel thanks to Lemma \ref{lem:lsym} (1) and 
\eqref{eq:l-sym}, and the term $(0,l,l)$ is
$0$. Therefore it is also true for $n$.
\end{proof}

Now we prove a part of Nekrasov's conjecture:
\begin{Theorem}\label{thm:regular}
We set
\begin{equation*}
\Fin(\ve_1,\ve_2,\vec{a};\q,\bbeta) \defeq
\ve_1 \ve_2 \log \Zin(\ve_1,\ve_2,\vec{a};\q,\bbeta).
\end{equation*}
Then 
$\Fin(\ve_1,\ve_2,\vec{a};\q,\bbeta)$ is regular at $(\ve_1,\ve_2)=(0,0)$.
\end{Theorem}

\begin{proof}
We omit $\bbeta$ by letting it be $1$ for brevity.

By \thmref{thm:key2}(1), we have
$\bZin_{0,d+1}-\bZin_{0,d}=0$ for $0 \leq d<r$.
We divide this equation by 
\(
   \Zin(\ve_1,\ve_2-\ve_1,\vec{a},t_1^{d-r/2} \q)
   \Zin(\ve_1-\ve_2,\ve_2,\vec{a},t_2^{d-r/2} \q).
\)
\begin{NB}
\begin{equation*}
\frac{\bZin_{0,d+1}-\bZin_{0,d}}
{\Zin(\ve_1,\ve_2-\ve_1,\vec{a},t_1^{d-r/2} \q)
\Zin(\ve_1-\ve_2,\ve_2,\vec{a},t_2^{d-r/2} \q)}
=0.
\end{equation*}
Then we have
\begin{equation*}
\begin{split}
  \sum_{\vec{k}}
 &  \frac{((t_1 t_2)^{d+1-r/2}\q)^{(\vec{k},\vec{k})/2}}
   {\prod_{\vec\alpha\in\Delta} l^{\vec{k}}_{\vec{\alpha}}(\ve_1,\ve_2,\vec{a})}
   \left(e^{(d+1)(\vec{k},\vec{a})}
   \frac{\Zin(\ve_1,\ve_2-\ve_1,\vec{a}+\ve_1\vec{k};t_1^{d+1-r/2}\q)
   \Zin(\ve_1-\ve_2,\ve_2,\vec{a}+\ve_2\vec{k};t_2^{d+1-r/2}\q)}
  {\Zin(\ve_1,\ve_2-\ve_1,\vec{a};t_1^{d-r/2} \q)
   \Zin(\ve_1-\ve_2,\ve_2,\vec{a};t_2^{d-r/2} \q)} \right.\\
& \quad - \left. e^{d(\vec{k},\vec{a})}
   \frac{\Zin(\ve_1,\ve_2-\ve_1,\vec{a}+\ve_1\vec{k};t_1^{d-r/2} \q)
   \Zin(\ve_1-\ve_2,\ve_2,\vec{a}+\ve_2\vec{k};t_2^{d-r/2}\q)}
  {\Zin(\ve_1,\ve_2-\ve_1,\vec{a};t_1^{d-r/2} \q)
   \Zin(\ve_1-\ve_2,\ve_2,\vec{a};t_2^{d-r/2} \q)}\right)=0.
\end{split}
\end{equation*}
\end{NB}

Let us expand $\Fin$ as
\begin{equation*}
   \Fin(\ve_1,\ve_2,\vec{a};\q) = \sum_n F_n(\ve_1,\ve_2,\vec{a}) \q^n.
\end{equation*}
We note that
\begin{equation*}
\begin{split}
&\frac{\Zin(\ve_1,\ve_2-\ve_1,\vec{a}+\ve_1\vec{k};t_1^{d+1-r/2}\q)
   \Zin(\ve_1-\ve_2,\ve_2,\vec{a}+\ve_2\vec{k};t_2^{d+1-r/2}\q)}
  {\Zin(\ve_1,\ve_2-\ve_1,\vec{a};t_1^{d-r/2} \q)
   \Zin(\ve_1-\ve_2,\ve_2,\vec{a};t_2^{d-r/2} \q)}\\
=\; &
\frac{\Zin(\ve_1,\ve_2-\ve_1,\vec{a}+\ve_1\vec{k};t_1^{d+1-r/2}\q)
   \Zin(\ve_1-\ve_2,\ve_2,\vec{a}+\ve_2\vec{k};t_2^{d+1-r/2}\q)}
{\Zin(\ve_1,\ve_2-\ve_1,\vec{a};t_1^{d+1-r/2}\q)
   \Zin(\ve_1-\ve_2,\ve_2,\vec{a};t_2^{d+1-r/2}\q)}\\
 & \quad \quad \quad \quad\quad\quad
   \times \frac{\Zin(\ve_1,\ve_2-\ve_1,\vec{a};t_1^{d+1-r/2}\q)
   \Zin(\ve_1-\ve_2,\ve_2,\vec{a};t_2^{d+1-r/2}\q)}
  {\Zin(\ve_1,\ve_2-\ve_1,\vec{a};t_1^{d-r/2} \q)
   \Zin(\ve_1-\ve_2,\ve_2,\vec{a};t_2^{d-r/2} \q)}\\
=\; &
\exp \left[\sum_{n \geq 1}\left( 
\frac{\Fin_n(\ve_1,\ve_2-\ve_1,\vec{a}+\ve_1\vec{k})-
\Fin_n(\ve_1,\ve_2-\ve_1,\vec{a})}{\ve_1(\ve_2-\ve_1)}
t_1^{(d+1-r/2)n}
\right.\right.\\
& \quad \quad \quad + \left.
\frac{\Fin_n(\ve_1-\ve_2,\ve_2,\vec{a}+\ve_1\vec{k})-
\Fin_n(\ve_1-\ve_2,\ve_2,\vec{a})}{\ve_2(\ve_2-\ve_1)}
t_2^{(d+1-r/2)n}
\right)\q^n \\
&\quad \quad \quad \quad +\left.\sum_{n \geq 1} 
\left(
\frac{\Fin_n(\ve_1,\ve_2-\ve_1,\vec{a})(t_1^n-1)t_1^{(d-r/2) n}}
{\ve_1(\ve_2-\ve_1)}-
\frac{\Fin_n(\ve_1-\ve_2,\ve_2,\vec{a})(t_2^n-1)t_2^{(d-r/2) n}}
{\ve_2(\ve_2-\ve_1)}
\right)\q^n \right].
\end{split}
\end{equation*}

By \eqref{eq:blow-up1} we get the following relations for $d=0,1$:
\begin{equation*}
\begin{split}
\sum_{n \geq 1}\left(
\frac{\Fin_n(\ve_1,\ve_2-\ve_1,\vec{a})(t_1^n-1)t_1^{-rn/2}}
{\ve_1(\ve_2-\ve_1)}-
\frac{\Fin_n(\ve_1-\ve_2,\ve_2,\vec{a})(t_2^n-1)t_2^{-rn/2}}
{\ve_2(\ve_2-\ve_1)}
\right)\q^n=\sum_n A_n \q^n\\
\sum_{n \geq 1}\left( 
\frac{\Fin_n(\ve_1,\ve_2-\ve_1,\vec{a})(t_1^n-1)t_1^{(1-r/2)n}}
{\ve_1(\ve_2-\ve_1)}-
\frac{\Fin_n(\ve_1-\ve_2,\ve_2,\vec{a})(t_2^n-1)t_2^{(1-r/2)n}}
{\ve_2(\ve_2-\ve_1)}
\right)\q^n=\sum_n B_n \q^n,
\end{split}
\end{equation*}
where the right hand sides comes from terms with $\vec{k} \ne 0$. In
particular, $A_n$ and $B_n$ are written by $F_m$ with $m<n$.
Solving the above, we get
\begin{equation*}
   F_n(\ve_1,\ve_2-\ve_1,\vec{a})=
   \frac{\ve_1(\ve_2-\ve_1)}{(t_1^n-1)(t_2^n-t_1^n)} t_1^{rn/2}
   (t_2^n A_n-B_n).
\end{equation*}
Hence if $F_m$, $m < n$ are regular at $(\ve_1,\ve_2)=(0,0)$, then
$F_n$ is also regular. As $F_0 = 0$, we get the assertion by the
induction on $n$.
\end{proof}

\subsection{The perturbative part}\label{subsec:pert}

We set
\begin{equation*}
\begin{split}
  & \gamma_{\ve_1,\ve_2}(x|\bbeta;\Lambda) \defeq
\frac{1}{2\ve_1\ve_2}\left(
-\frac{\bbeta}{6}\left(x+\frac{1}{2}(\ve_1+\ve_2)\right)^3 
+x^2\log(\bbeta\Lambda)\right)+\sum_{n \geq 1}\frac{1}{n}
\frac{e^{-\bbeta nx}}{(e^{\bbeta n\ve_1}-1)(e^{\bbeta n\ve_2}-1)},
\\
  & \widetilde{\gamma}_{\ve_1,\ve_2}(x|\bbeta;\Lambda)
\\
\defeq \; & 
  \gamma_{\ve_1,\ve_2}(x|\bbeta;\Lambda)
  + \frac{1}{\ve_1 \ve_2} \left(\frac{\pi^2 x}{6
  \bbeta}-\frac{\zeta(3)}{\bbeta^2} \right)
  +\frac{\ve_1+\ve_2}{2\ve_1 \ve_2} 
  \left( x \log (\bbeta \Lambda)+\frac{\pi^2}{6\bbeta} \right)+
  \frac{\ve_1^2+\ve_2^2+3\ve_1 \ve_2}{12 \ve_1 \ve_2} \log(\bbeta\Lambda)
\end{split}
\end{equation*}
for $\bbeta x>0$. Here $\Lambda = \q^{1/2r}$.
If we formally expand $\ve_1
\ve_2\widetilde{\gamma}_{\ve_1,\ve_2}(x|\bbeta;\Lambda)$ 
as a power series of $\ve_1,\ve_2$ around $\ve_1=\ve_2=0$, 
then each coefficient is a holomorphic function of
$x$. 
Indeed if we expand
\begin{equation*}
\frac{1}{(e^{\ve_1 \bbeta}-1)(e^{\ve_2 \bbeta}-1)}=
\sum_{m \geq 0} \frac{c_m}{m!}\bbeta^{m-2},
\end{equation*}
then
\begin{equation*}
\sum_{n \geq 1}\frac{1}{n}
\frac{e^{-\bbeta nx}}{(e^{\bbeta n\ve_1}-1)(e^{\bbeta n\ve_2}-1)}
=\sum_{m \geq 0} \frac{c_m}{m!}\bbeta^{m-2}
\mathrm{Li}_{3-m}(e^{-\bbeta x}),
\end{equation*}
where $\mathrm{Li}_{3-m}(e^{-\bbeta x})$ are polylogarithms.

\begin{NB}
We have
\begin{equation*}
\frac{1}{(e^{\bbeta \ve_1}-1)(e^{\bbeta \ve_2}-1)}
=\frac{1}{\ve_1 \ve_2}\sum_m \mathrm{Todd}_m({\Bbb C}^2)(-\bbeta)^{m-2}.
\end{equation*}
Therefore
\begin{equation*}
   \frac{c_m}{m!} = \frac{(-1)^m \Todd_m(\C^2)}{\ve_1\ve_2}.
\end{equation*}
Hence
\begin{equation*}
\begin{split}
&\widetilde{\gamma}_{\ve_1,\ve_2}(x|\bbeta;\Lambda)\\
=\; & 
\sum_{m\ge 0} \frac{(-1)^m \Todd_m(\C^2)}{\ve_1\ve_2} \bbeta^{m-2}
\mathrm{Li}_{3-m}(e^{-\bbeta x})
% \sum_{n \geq 1}\frac{1}{n}
% \frac{e^{-\bbeta nx}}{(e^{\bbeta n\ve_1}-1)(e^{\bbeta n\ve_2}-1)}
-\frac{\bbeta}{12 \ve_1\ve_2}\left(x-\frac{1}{2}K_{{\Bbb C}^2}\right)^3 \\
&+\frac{1}{\ve_1\ve_2}
\left[-\frac{1}{\bbeta^2}\zeta(3)+\left(x-\frac{K_{{\Bbb C}^2}}{2} \right)
\frac{\zeta(2)}{\bbeta}+
\left(\frac{x \cdot (x-K_{{\Bbb C}^2})}{2}+\mathrm{Todd}_2({\Bbb C}^2)\right)
\log(\bbeta\Lambda)\right].
\end{split}
\end{equation*}
Note that
\begin{equation*}
\sum_{n \geq 1}\frac{1}{n}
\frac{e^{-tnx}}{(e^{tn\ve_1}-1)(e^{tn\ve_2}-1)}
\end{equation*}
is holomorphic if $\bbeta\ve_1 \not \in {\Bbb R}\sqrt{-1}$,
$\bbeta\ve_2 \not \in {\Bbb R}\sqrt{-1}$ and $\bbeta x>0$.
Hence at least in this range
$\widetilde{\gamma}_{\ve_1,\ve_2}(x|\bbeta;\Lambda)$ is well defined.
On the other hand, I don't know the holomorphicity of
$\ve_1 \ve_2\widetilde{\gamma}_{\ve_1,\ve_2}(x|\bbeta;\Lambda)$ at
$\ve_1=\ve_2=0$.

\end{NB}
We extend the definition of $\gamma_{\ve_1,\ve_2}(x|\bbeta;\Lambda)$
to the range $\bbeta x<0$ by analytic continuation along a circle 
counter-clockwise way.
\begin{Proposition}\label{prop:perturbative}
If $\vec{a}$ is in a neighborhood of
the region $\bbeta\langle\vec{a},\vec{\alpha}\rangle>0$ for all $\vec{\alpha}
\in \Delta_+$, 
\begin{equation*}
\begin{split}
&\sum_{\vec{\alpha} \in \Delta}
\begin{aligned}[t]
\Bigl(
  & \widetilde{\gamma}_{\ve_1,\ve_2-\ve_1}
    (\langle\vec{a}+\vec{k}\ve_1,\vec{\alpha}\rangle|\bbeta;
    \Lambda e^{{(d-r/2)\bbeta\ve_1}/{2r}})
\\
  & \qquad\qquad +
  \widetilde{\gamma}_{\ve_1-\ve_2,\ve_2}
    (\langle \vec{a}+\vec{k}\ve_2,\vec{\alpha}\rangle|\bbeta;
    \Lambda e^{{(d-r/2)\bbeta\ve_2}/{2r}})
    -\widetilde{\gamma}_{\ve_1,\ve_2}
    (\langle\vec{a},\vec{\alpha}\rangle|\bbeta;\Lambda) \Bigr)
\end{aligned}
\\
= & 
\begin{aligned}[t]
  & - \frac{(4d - r)(r-1)}{48}(\ve_1+\ve_2)\bbeta
  - \frac{(\vec{k},\vec{k})}{2}\left( \log((\bbeta \Lambda)^{2r})
  - (d -\frac{r}2) (\ve_1+\ve_2)\bbeta\right)
\\
& \quad \quad\quad \quad\quad
  - d(\vec{k},\vec{a})\bbeta + \sum_{\vec{\alpha} \in \Delta}
\log (l_{\vec{\alpha}}^{\vec{k}}(\ve_1,\ve_2,\vec{a})).
\end{aligned}
\end{split}
\end{equation*}

\end{Proposition}

\begin{proof}
If $k=-l$ with $l>0$, then
\begin{equation}\label{eq:perturb-diff}
\begin{split}
&\frac{e^{-(x+k\ve_1)}}
{(e^{\ve_1}-1)(e^{\ve_2-\ve_1}-1)}+
\frac{e^{-(x+k\ve_2)}}
{(e^{\ve_2}-1)(e^{\ve_1-\ve_2}-1)}-
\frac{e^{-x}}
{(e^{\ve_1}-1)(e^{\ve_2}-1)}\\
=\; &e^{-x}\left(\frac{e^{\ve_1+\ve_2}
(e^{l\ve_1}-e^{l\ve_2})-(e^{(l+1)\ve_1}-e^{(l+1)\ve_2})
}
{(e^{\ve_1}-1)(e^{\ve_2}-1)(e^{\ve_2}-e^{\ve_1})}
-\frac{1}
{(e^{\ve_1}-1)(e^{\ve_2}-1)} \right)\\
=\; & e^{-x}\left(\frac{-\sum_{i=0}^{l-1}
e^{(i+1)\ve_1}e^{(l-i)\ve_2}+
\sum_{i=0}^{l} e^{i\ve_1}e^{(l-i)\ve_2}
}
{(e^{\ve_1}-1)(e^{\ve_2}-1)}
-\frac{1}
{(e^{\ve_1}-1)(e^{\ve_2}-1)} \right)\\
=\; &e^{-x}\frac{-\sum_{i=0}^{l-1}
e^{i\ve_1}e^{(l-i)\ve_2}+\sum_{i=0}^{l-1} e^{i\ve_1}}
{(e^{\ve_2}-1)}
=-e^{-x}\sum_{\substack{0 \leq i,j\\
i+j \leq l-1}}e^{i\ve_1}e^{j\ve_2}.
\end{split}
\end{equation}
Since 
\begin{equation*}
\frac{\ve_1 \ve_2 e^{-tx}}{(e^{t\ve_1}-1)(e^{t\ve_2}-1)}
\end{equation*}
is holomorphic at $(\ve_1,\ve_2)=(0,0)$,
\eqref{eq:perturb-diff} holds in 
${\Bbb C}[e^{-x}][[\ve_1,\ve_2]][\frac{1}{\ve_1\ve_2(\ve_1-\ve_2)}]$.

If $k>0$, then
\begin{equation*}
\begin{split}
&\frac{e^{-(x+k\ve_1)}}
{(e^{\ve_1}-1)(e^{\ve_2-\ve_1}-1)}+
\frac{e^{-(x+k\ve_2)}}
{(e^{\ve_2}-1)(e^{\ve_1-\ve_2}-1)}-
\frac{e^{-x}}
{(e^{\ve_1}-1)(e^{\ve_2}-1)}\\
=\;&e^{-x-\ve_1-\ve_2}
\left(\frac{
(e^{-k\ve_1}-e^{-k\ve_2})-
(e^{-(k-1)\ve_1}-e^{-(k-1)\ve_2})e^{-\ve_1-\ve_2}
}
{(e^{-\ve_1}-1)(e^{-\ve_2}-1)(e^{-\ve_1}-e^{-\ve_2})}
-\frac{1}
{(e^{-\ve_1}-1)(e^{-\ve_2}-1)} \right)\\
=\;& e^{-x}\left(\frac{\sum_{i=0}^{k-1}
e^{-i\ve_1}e^{-(k-1-i)\ve_2}-
\sum_{i=0}^{k-2} e^{-(i+1)\ve_1}e^{-(k-1-i)\ve_2}
}
{(e^{-\ve_1}-1)(e^{-\ve_2}-1)}
-\frac{1}
{(e^{-\ve_1}-1)(e^{-\ve_2}-1)} \right)\\
=\;&e^{-x}\frac{-\sum_{i=0}^{k-2}
e^{i\ve_1}e^{-(k-1-i)\ve_2}+\sum_{i=0}^{k-2} e^{-i\ve_1}}
{(e^{-\ve_2}-1)}
=-e^{-x}\sum_{\substack{0 \leq i,j\\
i+j \leq k-2}}e^{-(i+1)\ve_1}e^{-(j+1)\ve_2}.
\end{split}
\end{equation*}
Hence
\begin{equation}\label{eq:gamma}
\begin{split}
&\left.\sum_{n \geq 1}\frac{1}{n}
\frac{e^{-\bbeta n(\vec{a},\vec{\alpha})}}
  {(e^{\bbeta n\ve_1}-1)(e^{\bbeta n\ve_2}-1)}
\right|_{\substack{\vec{a} \to \vec{a}+\vec{k}\ve_1\\
\ve_1 \to \ve_1\\ \ve_2 \to \ve_2-\ve_1}}
+\left.\sum_{n \geq 1}\frac{1}{n}
\frac{e^{-\bbeta n(\vec{a},\vec{\alpha})}}
{(e^{\bbeta n\ve_1}-1)(e^{tn\ve_2}-1)}
\right|_{\substack{\vec{a} \to \vec{a}+\vec{k}\ve_2\\
\ve_1 \to \ve_1-\ve_2\\ \ve_2 \to \ve_2}}\\
&-\left.\sum_{n \geq 1}\frac{1}{n}
\frac{e^{-\bbeta n(\vec{a},\vec{\alpha})}}
{(e^{\bbeta n\ve_1}-1)(e^{\bbeta n\ve_2}-1)}
\right.
=\log(l_{\vec{\alpha}}^{\vec{k}}(\ve_1,\ve_2,\vec{a}))
\end{split}
\end{equation}
in ${\cal O}[[\ve_1,\ve_2]][\frac{1}{\ve_1\ve_2(\ve_1-\ve_2)}]$,
where ${\cal O}$ is the ring of holomorphic functions of $\bbeta\vec{a}$
in a neighborhood of $\bbeta\langle\vec{a},\vec{\alpha}\rangle>0$.
By the analytic continuation, 
\eqref{eq:gamma} also holds for the range
$\bbeta\langle\vec{a},\vec{\alpha}\rangle<0$. 
Combining the following equalities, we get our claim.
{\allowdisplaybreaks
\begin{gather}
\begin{aligned}[t]
&\frac{-(x+\ve_1 k+\frac{1}{2}\ve_2)^3}{12\ve_1(\ve_2-\ve_1)}+
\frac{-(x+\ve_2 k+\frac{1}{2}\ve_1)^3}{12\ve_2(\ve_1-\ve_2)}
\\
=\; & 
\frac{-(x+\frac{1}{2}(\ve_1+\ve_2))^3}{12\ve_1 \ve_2}
+ \frac{4k^2-4k+1}{16}x+\frac{8k^3-6k+2}{96}(\ve_1+\ve_2),
\end{aligned}
\\
\begin{aligned}[t]
&\frac{(x+k\ve_1)^2}{2\ve_1(\ve_2-\ve_1)}
 \log(\bbeta\Lambda e^{{(d-r/2)\bbeta\ve_1}/{2r}})
+\frac{(x+k\ve_2)^2}{2\ve_2(\ve_1-\ve_2)}
 \log(\bbeta\Lambda e^{{(d-r/2)\bbeta\ve_2}/{2r}})
\\
=\; & \frac{x^2}{2\ve_1\ve_2}\log(\bbeta \Lambda)
-\frac{k^2}{2}\log(\bbeta\Lambda)-
\frac{d-\frac{r}2}{2r} \left(\frac{k^2}{2}(\ve_1+\ve_2)+{kx} \right),
\end{aligned}
\\
  \frac{\pi^2(x+k\ve_1)}{6\bbeta\ve_1(\ve_2-\ve_1)}
  + \frac{\pi^2(x+k\ve_2)}{6\bbeta(\ve_1-\ve_2)\ve_2}
  = \frac{\pi^2 x}{6\bbeta\ve_1\ve_2},
\\
\begin{aligned}[t]
  & \frac{\ve_2(x+k\ve_1)\log(\bbeta\Lambda e^{{(d-r/2)\bbeta\ve_1}/{2r}})}
  {2\ve_1(\ve_2-\ve_1)}
  + \frac{\ve_1(x+k\ve_2)\log(\bbeta\Lambda e^{{(d-r/2)\bbeta\ve_2}/{2r}})}
  {2(\ve_1-\ve_2)\ve_2}
\\
  =\; & 
  \frac{(\ve_1+\ve_2)x\log(\bbeta\Lambda)}{2\ve_1\ve_2}
  + \frac{k}2 \log(\bbeta\Lambda)
  + \frac{(d - \frac{r}2)\bbeta x}{4r},
\end{aligned}
\\
\begin{aligned}[t]
   & \frac{-\ve_1^2+\ve_2^2+\ve_1\ve_2}
     {12\ve_1(\ve_2-\ve_1)}
   \log(\bbeta\Lambda e^{{(d-r/2)\bbeta\ve_1}/{2r}})
   + 
   \frac{\ve_1^2-\ve_2^2+\ve_1\ve_2}
     {12(\ve_1-\ve_2)\ve_2}
   \log(\bbeta\Lambda e^{{(d-r/2)\bbeta\ve_2}/{2r}})
\\
   =\; & \frac{\ve_1^2 + \ve_2^2 + 3\ve_1\ve_2}{\ve_1\ve_2}
   \log(\bbeta\Lambda) + \frac{(d - \frac{r}2)\bbeta(\ve_1+\ve_2)}{12r}.
\end{aligned}
\end{gather}
}
\end{proof}

\begin{NB}
Then we get
\begin{equation}
\begin{split}
&\sum_{\vec{\alpha} \in \Delta}
\left(
\gamma_{\ve_1,\ve_2-\ve_1}((\vec{a}+\vec{k}\ve_1,\vec{\alpha})|t;\Lambda
e^{\frac{dt}{2r}\ve_1})+
\gamma_{\ve_1-\ve_2,\ve_2}((\vec{a}+\vec{k}\ve_2,\vec{\alpha})|t;\Lambda 
e^{\frac{dt}{2r}\ve_2})
-\gamma_{\ve_1,\ve_2}((\vec{a},\vec{\alpha})|t;\Lambda) \right)\\
=& -\left( \frac{r}{2}(\vec{k},\vec{a})t+\frac{r(r-1)}{48}(\ve_1+\ve_2)t
+\frac{(\vec{k},\vec{k})}{2}\log((t \Lambda)^{2r})
+d \left(\frac{(\vec{k},\vec{k})}{2}(\ve_1+\ve_2)t+(\vec{k},\vec{a})t
\right) \right)\\
& \quad \quad\quad \quad+\sum_{\vec{\alpha} \in \Delta}
\log (l_{\vec{\alpha}}^{\vec{k}}(\ve_1,\ve_2,\vec{a})).
\end{split}
\end{equation}
\end{NB}

%Here I assume that $\log(-y)=-\pi \sqrt{-1}+\log y$.
%\begin{NB}
%The original formula is $\Li_3(-1/z)=\Li_3(-z)+...$.
%So we set $z=y e^{-\pi \sqrt{-1}}$.
%$\Li_k(z)=\sum_{n >0} z^k/n^k$, $k \geq 2$ 
%seems to be defined on $|z| \leq 1$.
%So we assume that $\Li_k$ is single valued on $|z| \leq 1$.
%Also I assume that $|y|=1$ ($\leftrightarrow$ $xt$ is pure imaginary).
%\end{NB}
\begin{NB}
Here I assume that $y=\bbeta x>0$. 
If we set $x=-\sqrt{-1}z, z>0$, then $y>0$ iff $\bbeta=s \sqrt{-1}$, $s>0$.
\end{NB} 
By using \eqref{formula:Li2} and \eqref{formula:Li3}, we see that
\begin{equation*}
\begin{split}
& \ve_1\ve_2 \left(\widetilde{\gamma}(x;\bbeta)
  + \widetilde{\gamma}(-x;\bbeta)\right)\\
=\; & 2\left( \frac{1}{\bbeta^2}(\Li_3(e^{-\bbeta x})-\zeta(3))+
\frac{x^2}{2}\log(\bbeta\Lambda)+\frac{\pi^2 x}{6\bbeta}\right)
-\frac{x^2 \pi \sqrt{-1}}{2}-\frac{\bbeta x^3}{6}
\\
& \quad +(\ve_1+\ve_2)\frac{x}{2}\pi \sqrt{-1}
-\frac{\ve_1^2+\ve_2^2+3\ve_1\ve_2}{6}
\left( 
\log \left(\frac{1-e^{-\bbeta x}}{\bbeta \Lambda} \right)
+ % \frac{\ve_1^2+\ve_2^2+3\ve_1\ve_2}{12}
  \frac{\bbeta x+\pi \sqrt{-1}}2\right)
  +\cdots.
\end{split}
\end{equation*}
By \subsecref{subsec:limit}, we get
\begin{equation*}
\begin{split}
& \ve_1\ve_2 \left(\widetilde{\gamma}(x;\bbeta)+
  \widetilde{\gamma}(-x;\bbeta)\right)\\ 
\underset{\bbeta \to 0}{\longrightarrow}&
\left(-\frac{1}{2}x^2 \log\left(\frac{x}{\Lambda} \right)+
\frac{3}{4}x^2 \right)
+\left(-\frac{1}{2}x^2 \log\left(\frac{-x}{\Lambda} \right)+
\frac{3}{4}x^2 \right)\\
& \quad +(\ve_1+\ve_2)\frac{x}{2}\pi \sqrt{-1}
-\frac{\ve_1^2+\ve_2^2+3\ve_1\ve_2}{6} 
\left(\log \left(\frac{x}{\Lambda} \right)
+ %\frac{\ve_1^2+\ve_2^2+3\ve_1\ve_2}{12}
  \frac{\pi \sqrt{-1}}2\right)+\cdots.
\end{split}
\end{equation*}
Therefore the perturbative part of the $K$-theoretic partition
function (see below) converges to that of the homological partition
function as $\bbeta\to 0$.

\subsection{The full partition functions}

By adding the perturbative term,
we define the {\it full partition functions\/} as
\begin{equation}
\begin{split}
Z(\ve_1,\ve_2,\vec{a};\q,\bbeta) \defeq & \exp 
\left(\sum_{\vec{\alpha} \in \Delta} 
-\widetilde{\gamma}_{\ve_1,\ve_2}(\langle\vec{a},\vec{\alpha}\rangle
|\bbeta;\Lambda) \right)
\Zin(\ve_1,\ve_2,\vec{a};\q,\bbeta)\\
\bZ_{k,d}(\ve_1,\ve_2,\vec{a};\q,\bbeta) \defeq & \exp 
\left(\sum_{\vec{\alpha} \in \Delta} 
-\widetilde{\gamma}_{\ve_1,\ve_2}(\langle\vec{a},\vec{\alpha}\rangle
|\bbeta;\Lambda) \right)
\Zin_{k,d}(\ve_1,\ve_2,\vec{a};\q,\bbeta).
\end{split}
\end{equation}

By \eqref{eq:blow-up1} and Proposition \ref{prop:perturbative},
we get
\begin{equation}\label{eq:blow-up4}
\begin{split}
 & \bZ_{k,d}(\ve_1,\ve_2,\vec{a};\q,\bbeta) =
  \sum_{\{\vec{k}\}=-{k}/{d}}
   \exp\left[ - \frac{(4d-r)(r-1)}{48} \bbeta(\ve_1+\ve_2) \right]\\ 
  &  \times 
  Z(\ve_1,\ve_2-\ve_1,\vec{a}+\ve_1\vec{k};t_1^{(d-\frac{r}{2})}\q,\bbeta)
   Z(\ve_1-\ve_2,\ve_2,\vec{a}+\ve_2\vec{k};t_2^{(d-\frac{r}{2})}\q,\bbeta).
\end{split}
\end{equation}
If we set $\ve_2=-\ve_1$, then
\begin{equation*}
%\begin{split}
 \bZ_{k,d}(\ve_1,-\ve_1,\vec{a};\q,\bbeta) =
  \sum_{\{\vec{k}\}=-{k}/{d}}
  Z(\ve_1,-2\ve_1,\vec{a}+\ve_1\vec{k};t_1^{(d-\frac{r}{2})}\q,\bbeta)
   Z(2\ve_1,-\ve_1,\vec{a}-\ve_1\vec{k};t_1^{-(d-\frac{r}{2})}\q,\bbeta).
%\end{split}
\end{equation*}

We expand $F$ as
\begin{equation*}
\begin{split}
&F(\ve_1,\ve_2,\vec{a};\q,\bbeta)
  = \ve_1 \ve_2 \log Z(\ve_1,\ve_2,\vec{a};\q,\bbeta)\\
=\;&F_0(\vec{a};\q,\bbeta)+(\ve_1+\ve_2)H(\vec{a};\q,\bbeta)+
(\ve_1+\ve_2)^2 G(\vec{a};\q,\bbeta)+\ve_1 \ve_2
F_1(\vec{a};\q,\bbeta)+\cdots.
\end{split}
\end{equation*}
By Lemma \ref{lem:blow-up}, $H(\vec{a};\q,\bbeta)$ 
comes from the perturbative part:
\begin{equation*}
\begin{split}
&H(\vec{a};\q,\bbeta)=\sum_{\alpha \ne \beta}
\left( \frac{1}{2\bbeta} \Li_2(e^{-\bbeta(a_\alpha-a_\beta)})+
\frac{\bbeta}{8}(a_\alpha-a_\beta)^2 -
\frac{1}{2}(a_\alpha-a_\beta)\log(\bbeta
\Lambda)-\frac{\pi^2}{12\bbeta}\right)\\ 
=& 
-\sum_{\alpha<\beta} \pi \sqrt{-1}\frac{(a_\alpha-a_\beta)}{2}=
-\pi \sqrt{-1} \langle \vec{a},\rho \rangle.
\end{split}
\end{equation*}
\begin{NB}
Here I assume that $x\bbeta>0$. I changed the sign. (Jan. 18). 
\end{NB}

Let us derive equations for $F_0$. We use
{\allowdisplaybreaks
\begin{equation*}
\begin{split}
   & \frac{F_0(\vec{a}+\ve_1\vec{k}; t_1^{\left(d-r/2\right)}\q,\bbeta)}
   {\ve_1(\ve_2-\ve_1)}
    +\frac{F_0(\vec{a}+\ve_2\vec{k}; t_2^{\left(d-r/2\right)}\q,\bbeta)}
    {(\ve_1-\ve_2)\ve_2}
\\
%  =\; &
  & \qquad\qquad\qquad =
  \frac{1}{\ve_1\ve_2}F_0
  - \left[%\sum_{p,q}
    \frac{\partial^2 F_0
    }{(\partial\log \q)^2} \frac{\bbeta^2}2\left(d-\frac{r}2\right)^2
    +
    %\sum_p
    \frac{\partial^2 F_0
    }{\partial\log \q\partial a^l}\bbeta\left(d-\frac{r}2\right) k^l
    +
    \frac{\partial^2 F_0
    }{\partial a^l\partial a^m}
      \frac{k^l k^m}{2}
  \right]
  + \cdots,
% \\
%  & -(\ve_1+\ve_2)
%  \left[\frac{\partial^3 F_0
%    }{\partial\tau_p\partial\tau_q\partial\tau_r}
%    \frac{t_p t_q t_r}{3!}
%    +
%   \frac{\partial^3 F_0
%   }{\partial\tau_p\partial\tau_q \partial a^l}
%   \frac{k^l t_p t_q}{2}
%   +
%   \frac{\partial^3 F_0
%   }{\partial \tau_p \partial a^l\partial a^m}\frac{k^l k^m t_p}{2}
%   +
%   \frac{\partial^3 F_0
%   }{\partial a^l\partial a^m\partial a^n}
%   \frac{k^l k^m k^n}{3!}
%   \right]
% \\
%   & -\left((\ve_1 + \ve_2)^2 - \ve_1\ve_2\right)
%   \Biggl[
%   \begin{aligned}[t]
%   & \frac{\partial^4 F_0
%   }{\partial\tau_p\partial\tau_q\partial\tau_r\partial\tau_s}
%   \frac{t_pt_qt_rt_s}{4!}+
%   \frac{\partial^4 \Fz%(\vec{a};\q)
%   }{\partial\tau_p\partial\tau_q\partial\tau_r\partial a^l}
%   \frac{k^l t_pt_qt_r}{3!}+
%   \frac{\partial^4 \Fz%(\vec{a};\q)
%   }{\partial\tau_p\partial\tau_q \partial a^l\partial a^m}
%   \frac{k^l k^m t_pt_q}{2\cdot 2}
% \\
%   & \qquad
%   + \frac{\partial^4 \Fz%(\vec{a};\q)
%   }{\partial\tau_p\partial a^l
%   \partial a^m \partial a^n}\frac{k^l k^m k^n t_p}{3!}+
%   \frac{\partial^4 \Fz%(\vec{a};\q)
%   }{\partial a^l\partial a^m \partial a^n \partial a^o}
%   \frac{k^l k^m k^n k^o}{4!}
%   \Biggr] +\cdots
%  \end{aligned}
\\
   & \frac{\ve_2 H(\vec{a}+\ve_1\vec{k}; t_1^{\left(d-r/2\right)}\q,\bbeta)}
   {\ve_1(\ve_2-\ve_1)}
    +\frac{\ve_1 H(\vec{a}+\ve_2\vec{k};t_2^{\left(d-r/2\right)}\q,\bbeta)}
   {(\ve_1-\ve_2)\ve_2}
\\
%  =
%  \; &
  & \qquad\qquad =
  \frac{\ve_1+\ve_2}{\ve_1\ve_2} H%(\vec{a};\q,\vec{\tau})
  + \left[\frac{\partial H
    }{\partial\log \q} \bbeta\left(d-\frac{r}2\right)
    +
    \frac{\partial H
    }{\partial a^l} k^l
  \right]
  + \cdots
  = \frac{\ve_1+\ve_2}{\ve_1\ve_2} H%(\vec{a};\q,\vec{\tau})
  - \pi\sqrt{-1}\langle\vec{k},\rho\rangle  + \cdots
  ,
% \\
%  & -\ve_1\ve_2
%  \left[\frac{\partial^3 H%(\vec{a};\q)
%    }{\partial\tau_p\partial\tau_q\partial\tau_r}
%    \frac{t_pt_qt_r}{3!}
%    +
%   \frac{\partial^3 H%(\vec{a};\q)
%   }{\partial\tau_p\partial\tau_q\partial a^l}
%   \frac{k^l t_pt_q}{2}
%   +
%   \frac{\partial^3 H%(\vec{a};\q)
%   }{\partial\tau_p\partial a^l\partial a^m}\frac{k^l k^m t_p}{2}
%   +
%   \frac{\partial^3 H%(\vec{a};\q)
%   }{\partial a^l\partial a^m\partial a^n}
%   \frac{k^l k^m k^n}{3!}
%   \right] + \cdots
\\
   & \frac{\ve_2^2 G(\vec{a}+\ve_1\vec{k};t_1^{\left(d-\frac1r\right)}\q,\bbeta)}
  {\ve_1(\ve_2-\ve_1)}
    +\frac{\ve_1^2 G(\vec{a}+\ve_2\vec{k};t_2^{\left(d-\frac1r\right)}\q,\bbeta)}
  {(\ve_1-\ve_2)\ve_2}
%\\
  =
%  \; &
  \frac{(\ve_1+\ve_2)^2-\ve_1\ve_2}{\ve_1\ve_2} G%(\vec{a};\q)
%   + (\ve_1+\ve_2)
%   \left[\frac{\partial G%(\vec{a};\q)
%     }{\partial\log \q} t{\left(d-\frac1r\right)}
%     +
%     \frac{\partial G%(\vec{a};\q)
%     }{\partial a^l} k^l \right]
% \\
%  & + \ve_1\ve_2
%  \left[\frac{\partial^2 G%(\vec{a};\q)
%    }{\partial\tau_p\partial\tau_q}
%    \frac{t_pt_q}{2}
%    +
%   \frac{\partial^2 G%(\vec{a};\q)
%   }{\partial\tau_p\partial a^l}{k^l t_p}
%   +
%   \frac{\partial^2 G%(\vec{a};\q)
%   }{\partial a^l\partial a^m}\frac{k^l k^m}{2}
%   \right]
  + \cdots,
\\
   & F_1(\vec{a}+\ve_1\vec{k}; t_1^{\left(d-\frac1r\right)}\q,\bbeta)
    + F_1(\vec{a}+\ve_2\vec{k}; t_2^{\left(d-\frac1r\right)}\q,\bbeta)
% \\
   = %\; &
  2 F_1
%   + (\ve_1+\ve_2)
%   \left[\frac{\partial F_1
%     }{\partial\tau_p} t_p
%     +
%     \frac{\partial F_1
%     }{\partial a^l} k^l
%   \right]
% \\
%  & + \left((\ve_1+\ve_2)^2 - 2\ve_1\ve_2\right)
%  \left[\frac{\partial^2 F_1
%    }{\partial\tau_p\partial\tau_q}
%    \frac{t_pt_q}{2}
%    +
%   \frac{\partial^2 F_1
%   }{\partial\tau_p\partial a^l}
%   {k^l t_p}
%   +
%   \frac{\partial^2 F_1
%   }{\partial a^l\partial a^m}\frac{k^l k^m}{2}
%   \right]
  + \cdots,
\end{split}
\end{equation*}
where $F_0$, $H$, $G$, $F_1$ and their derivatives are evaluated at
$\vec{a}, \q$ in the right hand sides.

By \thmref{thm:key2}(1) and \eqref{eq:blow-up4}, we have
\begin{equation}\label{eq:genus0}
\begin{split}
&\exp(F_1)=
\exp(2F_1-G) \sum_{\{\vec{k}\}=0} (-1)^{-\langle \vec{k},\rho \rangle} \times\\
&\quad \quad \exp \left[-\frac{\partial^2 F_0}{(\partial \log \q)^2}
\frac{\bbeta^2}{2}\left(d-\frac{r}{2} \right)^2-
\frac{\partial^2 F_0}{\partial \log \q \partial a^l}
\bbeta \left(d-\frac{r}{2} \right)k^l-
\frac{\partial^2 F_0}{\partial a^l \partial a^m}\frac{k^l k^m}{2}
\right]
\end{split}
\end{equation}
for $0 \leq d \leq r$.
\begin{NB}
If $r$ is even, then by setting $d=r/2$, we have
\begin{equation}
\exp(G-F_1)=\sum_{\vec{k}} (-1)^{-\langle \vec{k},\rho \rangle}
\exp \left[-\sum_{l,m}
\frac{\partial^2 F_0}{\partial a^l \partial a^m}\frac{k^l k^m}{2}
\right]
\end{equation}
\end{NB}

\begin{NB}
Since $c_1=0$, 
$(-1)^{-\langle \vec{k},\rho \rangle}=
(-1)^{\langle \vec{k},\rho \rangle}$.
If $c_1=k$, then $(-1)^{-\langle \vec{k},\rho \rangle}=
(-1)^{\langle \vec{k},\rho \rangle}(-1)^{k(r-k)}$. 
\end{NB}

We set
\begin{equation*}
\begin{split}
\tau_{lm}(\bbeta) \defeq &-\frac{1}{2\pi \sqrt{-1}}
\frac{\partial^2 F_0(\vec{a};\q,\bbeta)}{\partial a^l \partial a^m}
% ,
% \\
% u_2(\bbeta)
% \defeq & -\frac{\partial F_0(\vec{a};\q,\bbeta)}{\partial \log \q}
% =-\frac{1}{2r}\Lambda\frac{\partial F_0(\vec{a};\q,\bbeta)}{\partial
% \Lambda }
.
\end{split}
\end{equation*}
\begin{NB}
I have slightly change the contact term equation. I also drop the
definition of $u_2$. (May 20)
\end{NB}
Then \eqref{eq:genus0} can be written as
\begin{equation}\label{eq:Kcontact}
\begin{split}
& \exp(G - F_1) = 
  \exp\left[- \frac{\partial^2 F_0}{(\partial \log \q)^2}
    \frac{\bbeta^2}{2}\left(d-\frac{r}{2} \right)^2
  \right]
  \Theta_E \left(\left. - \frac1{2\pi\sqrt{-1}}\frac{\partial^2 F_0}
      {\partial\log\q\partial\vec{a}}
  \bbeta \left(d-\frac{r}{2} \right)
  \right| \tau(\bbeta) \right)
\end{split}
\end{equation}
where $\Theta_E$ is the Riemann theta function with the characteristic
${}^t\left(\frac{1}{2},\frac{1}{2},\dots,\frac{1}{2} \right)$. 
(See \cite[Appendix~B]{lecture} for convention.)

As the left hand side is independent of $d$, we have
\begin{equation}\label{eq:Kcontact2}
\begin{split}
  & \exp\left[- \frac{\partial^2 F_0}{(\partial \log \q)^2}
    \frac{\bbeta^2}{2}\left(d-\frac{r}{2} \right)^2
  \right]
  \Theta_E \left(\left. - \frac1{2\pi\sqrt{-1}}\frac{\partial^2 F_0}
      {\partial\log\q\partial\vec{a}}
  \bbeta \left(d-\frac{r}{2} \right)
  \right| \tau(\bbeta) \right)
\\
 = \; & \exp\left[- \frac{\partial^2 F_0}{(\partial \log \q)^2}
    \frac{\bbeta^2}{2}\left(d-1-\frac{r}{2} \right)^2
  \right]
  \Theta_E \left(\left. - \frac1{2\pi\sqrt{-1}}\frac{\partial^2 F_0}
      {\partial\log\q\partial\vec{a}}
  \bbeta \left(d-1-\frac{r}{2} \right)
  \right| \tau(\bbeta) \right)
% & \exp \left[-\frac{\partial^2 F_0}{(\partial \log \q)^2}
% \frac{\bbeta^2}{2}\left(d-1-\frac{r}{2} \right)^2 \right] \times\\
% &\quad \quad
% \sum_{\{\vec{k}\}=0} (-1)^{-\langle \vec{k},\rho \rangle}
%  \exp \left[
% - \frac{\partial^2 F_0}{\partial \log \q \partial a^l}
% \bbeta \left(d-1-\frac{r}{2} \right)k^l-
% \frac{\partial^2 F_0}{\partial a^l \partial a^m}\frac{k^l k^m}{2}
% \right]\\
% =& \exp \left[-\frac{\partial^2 F_0}{(\partial \log \q)^2}
% \frac{\bbeta^2}{2}\left(d-\frac{r}{2} \right)^2 \right] \times\\
% &\quad \quad
% \sum_{\{\vec{k}\}=0} (-1)^{-\langle \vec{k},\rho \rangle}
%  \exp \left[
% -\frac{\partial^2 F_0}{\partial \log \q \partial a^l}
% \bbeta \left(d-\frac{r}{2} \right)k^l-
% \frac{\partial^2 F_0}{\partial a^l \partial a^m}\frac{k^l k^m}{2}
% \right]
\end{split}
\end{equation}
for $1 \leq d \leq r$.
If $d \not= {(r+1)}/{2}$, then
this is a differential equation for $F_0$ whose 
solution is determined by the perturbative part, as in
\corref{cor:recursive}.
\begin{NB}
If $d={(r+1)}/{2}$, then the equation is trivial.
\end{NB}
 
By taking the difference with respect to $d$ again,
we recover the contact term equation for the homological partition
function (see \cite[\S7]{part1} and \cite[\S5.3]{lecture}).
\begin{NB}
\begin{equation*}
\bbeta^2 \sum_{\{\vec{k}\}=0} (-1)^{\langle \vec{k},\rho \rangle}
\left[\frac{\partial^2 F_0}{(\partial \log \q)^2}-\left(
\frac{\partial^2 F_0}{\partial \log \q \partial a^l} k^l \right)^2
\right]
\exp \left(-
\frac{\partial^2 F_0}{\partial a^l \partial a^m}\frac{k^l k^m}{2} 
\right)
=O(\bbeta^3).
\end{equation*}
\end{NB}

\section{Genus $1$ part}\label{sec:genus1}

The coefficients $G$, $F_1$ are called {\it genus $1$ part\/} of the
partition function in the physics literature. It is known that $F_1$
is Gromov-Witten invariants for certain noncompact Calabi-Yau
$3$-folds. We determine these in terms of $\tau$ in this section (see
\eqref{eq:genus1}).

\subsection{Coefficients of $\Lambda^k$} We first prepare a result used
in the next subsection.

We set $\zeta_{\alpha,\beta}\defeq\frac{1}{1-e^{-(a_\alpha-a_\beta)\bbeta}}$.
Then we see that
\begin{equation}
\begin{split}
\frac{1}{1-e^{(i\ve_1+j\ve_2)\bbeta}e^{-(a_\alpha-a_\beta)\bbeta}}&=
\frac{1}{1-e^{-(a_\alpha-a_\beta)\bbeta}+e^{-(a_\alpha-a_\beta)\bbeta}
(1-e^{(i\ve_1+j\ve_2)\bbeta})}\\
&=\frac{1}{1-e^{-(a_\alpha-a_\beta)\bbeta}}\frac{1}{1-\frac{-e^{-(a_\alpha-a_\beta)\bbeta}}
{1-e^{-(a_\alpha-a_\beta)\bbeta}}(1-e^{(i\ve_1+j\ve_2)\bbeta})}\\
&=\zeta_{\alpha,\beta}\frac{1}{1-(1-\zeta_{\alpha,\beta})
(1-e^{(i\ve_1+j\ve_2)\bbeta})}\\
&=\zeta_{\alpha,\beta}\sum_{k \geq 0}(1-\zeta_{\alpha,\beta})^k
(1-e^{(i\ve_1+j\ve_2)\bbeta})^k
\end{split}
\end{equation}
and
\begin{equation}
\frac{1}{1-e^{(i\ve_1+j\ve_2)\bbeta}e^{-(a_\beta-a_\alpha)\bbeta}}
=(1-\zeta_{\alpha,\beta})\sum_{k \geq 0}\zeta_{\alpha,\beta}^k
(1-e^{(i\ve_1+j\ve_2)\bbeta})^k.
\end{equation}
Hence 
$$
\Zin \in \left({\Bbb C}[\zeta_{\alpha,\beta}][[\ve_1,\ve_2]] \otimes 
{\Bbb C}(\ve_1,\ve_2) \right)
[[\prod_{\alpha<\beta}\zeta_{\alpha,\beta}\Lambda]]
$$
In particular,
\begin{equation}\label{eq:zeta1}
\Fin \in {\Bbb C}[\zeta_{\alpha,\beta}]
[[\prod_{\alpha<\beta}\zeta_{\alpha,\beta}\Lambda]][[\ve_1,\ve_2]].
\end{equation}
Since
\begin{equation}
\begin{split}
\pd{\zeta_{\alpha,\beta}^n}{a_l}=& n \zeta_{\alpha,\beta}^{n-1}\pd{\zeta_{\alpha,\beta}}{a_l}\\
=\; & n \zeta_{\alpha,\beta}^{n-1}
\frac{-e^{-(a_\alpha-a_\beta)\bbeta}}{(1-e^{-(a_\alpha-a_\beta)\bbeta})^2}
\pd{(a_\alpha-a_\beta)\bbeta}{a_l}\\
= & \bbeta n\zeta_{\alpha,\beta}^n(1-\zeta_{\alpha,\beta})(e_\alpha-e_\beta,e_l),
\end{split}
\end{equation}
we also have
\begin{equation}\label{eq:zeta2}
\pd{}{a_l} \pd{\Fin_0}{ a_m} \in {\Bbb C}[\zeta_{\alpha,\beta}][[\prod_{\alpha<\beta}\zeta_{\alpha,\beta}\Lambda]].
\end{equation}

\subsection{Genus $1$ parts as modular forms}

Assume that $0<k<r$.
By \thmref{thm:key2}(2), we get
\begin{equation*}
\sum_{\{\vec{k}\}=-{k}/{r}}
Z(\ve_1,\ve_2-\ve_1,\vec{a}+\ve_1\vec{k};t_1^{d-\frac{r}{2}}\q,\bbeta)
Z(\ve_1-\ve_2,\ve_2,\vec{a}+\ve_2\vec{k};t_2^{d-\frac{r}{2}}\q,\bbeta)=0
\end{equation*}
for $0<d<r$.

As in the derivation of \eqref{eq:genus0} we have
{\allowdisplaybreaks
\begin{equation*}
\begin{split}
  0 &= \sum_{\{\vec{k}\} = -\frac{k}r}\!\!
  \begin{aligned}[t]
  & 
  \exp\Biggl[
  - \frac{\partial^2 F_0}{\partial a^l\partial a^m}
  \frac{k^l k^m}2-
\frac{\partial^2 F_0}{\partial \log \q \partial a^l}
\bbeta \left(d-\frac{r}{2} \right)k^l-\frac{\partial^2 F_0}{(\partial \log \q)^2}
\frac{\bbeta^2}{2}\left(d-\frac{r}{2} \right)^2
  + \frac{\partial H}{\partial a^l}k^l
    \\
  & \qquad + 
  (\ve_1+\ve_2)\left\{
    -\frac{\partial^3 F_0}{\partial a^l\partial a^m\partial a^n}
  \frac{k^l k^m k^n}{3!}\right.-\frac{\partial^3 F_0}{\partial a^l\partial a^m\partial \log \q}
  \frac{k^l k^m}{2!}\bbeta\left(d-\frac{r}{2} \right)\\
& \qquad -
\frac{\partial^3 F_0}{\partial a^l\partial(\log \q)^2}
  \frac{k^l}{2!}\bbeta^2\left(d-\frac{r}{2} \right)^2-
\frac{\partial^3 F_0}{\partial(\log \q)^3}\frac{1}{3!}
 \bbeta^3\left(d-\frac{r}{2} \right)^3\\
& \qquad  + \left.
    \frac{\partial (G+ F_1)}{\partial a^l} k^l
+\frac{\partial (G+ F_1)}{\partial \log \q}\bbeta\left(d-\frac{r}{2} \right)
  \right\} + \cdots
  \Biggr].
  \end{aligned}
\end{split}
\end{equation*}
Hence we have
\begin{equation*}
\begin{split}
  0 &= 
 \sum_{\{\vec{k}\} = -\frac{k}r}\!\!
  (-1)^{-\langle \vec{k},\vec{\rho} \rangle}
  \begin{aligned}[t]
  & 
  \exp\Biggl[
  - \frac{\partial^2 F_0}{\partial a^l\partial a^m}
  \frac{k^l k^m}2-
\frac{\partial^2 F_0}{\partial \log \q \partial a^l}
\bbeta \left(d-\frac{r}{2} \right)k^l
    \\
  & \qquad + 
  (\ve_1+\ve_2)\left\{
    -\frac{\partial^3 F_0}{\partial a^l\partial a^m\partial a^n}
  \frac{k^l k^m k^n}{3!}\right.-\frac{\partial^3 F_0}{\partial a^l\partial a^m\partial \log \q}
  \frac{k^l k^m}{2!}\bbeta\left(d-\frac{r}{2} \right)\\
& \qquad -
\frac{\partial^3 F_0}{\partial a^l\partial(\log \q)^2}
  \frac{k^l}{2!}\bbeta^2\left(d-\frac{r}{2} \right)^2+ \left.
    \frac{\partial (G+ F_1)}{\partial a^l} k^l
  \right\} + \cdots
  \Biggr].
  \end{aligned}
\end{split}
\end{equation*}
Setting $\ve_1=\ve_2=0$, we have
\begin{equation*}
  0 = 
 \sum_{\{\vec{k}\} = -\frac{k}r}\!\!
  (-1)^{-\langle \vec{k},\vec{\rho} \rangle}
  \exp\Biggl[
  - \frac{\partial^2 F_0}{\partial a^l\partial a^m}
  \frac{k^l k^m}2-
\frac{\partial^2 F_0}{\partial \log \q \partial a^l}
\bbeta \left(d-\frac{r}{2} \right)k^l\Biggr].
\end{equation*}

By looking at the coefficient of $\ve_1+\ve_2$, we get
\begin{equation*}
\begin{split}
  0 &= 
 \sum_{\{\vec{k}\} = -\frac{k}r}\!\!
  (-1)^{-\langle \vec{k},\vec{\rho} \rangle}
  \begin{aligned}[t]
  & 
  \exp\Biggl[
  - \frac{\partial^2 F_0}{\partial a^l\partial a^m}
  \frac{k^l k^m}2-
\frac{\partial^2 F_0}{\partial \log \q \partial a^l}
\bbeta \left(d-\frac{r}{2} \right)k^l \Biggr] 
    \\
  & \qquad \times 
  \left\{
    -\frac{\partial^3 F_0}{\partial a^l\partial a^m\partial a^n}
  \frac{k^l k^m k^n}{3!}\right.-\frac{\partial^3 F_0}{\partial a^l\partial a^m\partial \log \q}
  \frac{k^l k^m}{2!}\bbeta\left(d-\frac{r}{2} \right)\\
& \qquad -
\frac{\partial^3 F_0}{\partial a^l\partial(\log \q)^2}
  \frac{k^l}{2!}\bbeta^2\left(d-\frac{r}{2} \right)^2+ \left.
    \frac{\partial (G+ F_1)}{\partial a^l} k^l
  \right\}.
  \end{aligned}
\end{split}
\end{equation*}

 From now on, we assume that $r=2$ and $d=1$.
We set $\vec{a}:=(-a,a)$.
Then

\begin{equation*}
\sum_{\{\vec{k}\}=-\frac{1}{2}}(-1)^{-\langle\vec{k},\rho\rangle}
e^{\pi\sqrt{-1}\tau k^2}
\left(-\frac{\partial^3 F_0}{\partial a^3}\frac{k^3}{3!}
+\frac{\partial (F_1+G)}{\partial a}k \right)=0.
\end{equation*}
This equation can be rewritten as
\begin{equation*}
\frac{\partial }{\partial a}(G+F_1)=-\frac{1}{3}\frac{\partial }{ \partial a}
\log \left(\pd{}{\xi}\theta_{11}(\xi|\tau)_{|\xi=0} \right).
\end{equation*}
By Jacobi's derivative formula, we get
\begin{equation*}
G+F_1=-\frac{1}{3} \log \left(
-2\pi q^{\frac{1}{4}}\prod_{d=1}^{\infty}(1-q^{2d})^3 \right)+C
\end{equation*}
where $C$ is a function on $\Lambda$.
Combining this with \eqref{eq:genus0} (with $r=2$, $d=1$), we see that
\begin{equation*}
\exp(F_1)=\frac{C'}{\eta(\tau/2)}
\end{equation*}
where $C'$ is a function on $\Lambda$.
On the other hand,
\begin{equation*}
\exp(F_1)\eta(\tau/2) \in {\Bbb C}[\zeta_{1,2}][[\zeta_{1,2}\Lambda]].
\end{equation*}

\begin{proof}
We note the following relations:
\begin{equation*}
F_1=\frac{1}{12}\log \left(\frac{1-e^{-a\bbeta}}{\bbeta\Lambda} \right)+
\frac{1}{12}\log \left(\frac{1-e^{a\bbeta}}{\bbeta\Lambda} \right)+
\Fin_1,
\end{equation*}

\begin{equation*}
\begin{split}
\log \eta(\tau/2)^{-1}&=-\frac{1}{24}\pi\sqrt{-1}\tau +O(q),\\
\tau&=\frac{\sqrt{-1}}{\pi}2 \left(
\log \left(\frac{1-e^{-a\bbeta}}{\bbeta\Lambda} \right)+
\log\left(\frac{1-e^{a\bbeta}}{\bbeta\Lambda} \right) \right)-\frac{1}{2\pi \sqrt{-1}}
\frac{\partial^2 \Fin_0}{\partial a^2}.
\end{split}
\end{equation*}
Since $\Fin_1$, $\partial^2 \Fin_0/\partial a^2 \in
{\Bbb C}[\zeta_{1,2}][[\zeta_{1,2}\Lambda]]$ 
(\eqref{eq:zeta1},\eqref{eq:zeta2}),
we get our claim.
\end{proof}

Hence $C'$ is a constant. By the same proof, we also see that $C'=1$.   
Therefore we get the following equalities:
\begin{equation}\label{eq:genus1}
\begin{split}
\exp(F_1) & =\frac{1}{\eta(\tau/2)},\\
\exp(G) &=q^{-1/24} \prod_{d=1}^{\infty}(1-q^{2d-1}).
\end{split}
\end{equation}

\appendix
\section{Normal varieties}

\begin{Lemma}\label{lem:morphism}
Let $X$ be a normal variety over ${\Bbb C}$.
Let $f:X \to Y$ be a proper and surjective morphism
to a normal variety $Y$
and $g:X \to Z$ a morphism
to a variety $Z$.
Assume that $g(f^{-1}(y))$ is a point for every $y \in Y$.
Then we have a unique morphism $h:Y \to Z$ such that
$h \circ f=g$.
\end{Lemma}

\begin{proof}
We note that $f$ is factorized
$X \hookrightarrow X \times Z \overset{f \times 1_Z}{\to} Y \times Z \to Y$. 
Since $f \times 1_Z$ is proper,
$W:=(f \times 1_Z)(X)$ is a closed subvariety of
$Y \times Z$. Since $X \to W$ is surjective and
$X \to W \to Y$ is proper, $W \to Y$ is also proper.
Since $g(f^{-1}(y))$ is a point for every $y \in Y$,
$W \to Y$ is injective.
Hence it is finite and birational.
Since $Y$ is normal, we conclude that $W \cong Y$.
Therefore we have a desired morphism $h:Y \to Z$.
\end{proof}

\begin{Corollary}\label{cor:equivalence}
Let $X$ be a normal variety and 
$R$ an equivalence relation on $X({\Bbb C})$
(the set of ${\Bbb C}$-valued points of $X$).
Assume that there is a proper map
$g:X \to Y$ such that $Y({\Bbb C})=X({\Bbb C})/R$.
Then there is a proper morphism
$f:X \to M$ such that
$M$ is normal, $M({\Bbb C})=X({\Bbb C})/R$
and there is a unique morphism
$h:M \to Y$ with $g=h \circ f$.
We call $M$ the quotient of $X$ by $R$.
\end{Corollary}

\section{Polylogarithms}

\subsection{}
We define the $k$th polylogarithm inductively by
\begin{equation*}
\begin{gathered}
\Li_k(0) = 0, \qquad \Li_0(w) = \frac{w}{1-w},
\\
  \deriv{w} \Li_k(w) % = \sum_{n=1}^\infty \frac{w^{n-1}}{n^{k-1}}
  = \frac1w \Li_{k-1}(w),
\end{gathered}
\end{equation*}
where $w < 1$. We have
\begin{equation*}
   \Li_1(w) = \int \frac1w \Li_{0}(w) dw
   = - \log (1 - w).
\end{equation*}
We also have
\begin{equation*}
\begin{gathered}
\Li_{-1}(w) = w \deriv{w}\Li_0(w) = \frac{w}{(1-w)^2},
\quad
   \Li_{-2}(w) = w \deriv{w}\Li_{-1}(w) = \frac{w + w^2}{(1-w)^3},
\\
   \Li_{-k}(w) = \frac{w P_{k}(w)}
   {(1-w)^{k+1}} \qquad (k\in\Z_{\ge 0}),
\end{gathered}
\end{equation*}
where $P_0(w) = 1$ and $P_{k}(w)$ is a polynomial of degree $k-1$ for
$k > 0$, satisfying the recursive system
\begin{equation*}
   P_{k+1}(w) = (1 + kw) P_{k}(w) + w(1-w) P_{k}'(w).
\end{equation*}
In particular, $\Li_{-k}(w)$ is defined for $w\neq 1$.
For $|w|<1$, $\Li_k$ can be expressed by a series
\begin{equation*}
  \Li_k(w) = \sum_{n=1}^\infty \frac{w^n}{n^k}.
\end{equation*}

\subsection{Inversion formulas}
We need to relate $\Li_k(1/w)$ to $\Li_k(w)$. For $k = -l$ with
$l > 0$, we have
\begin{equation*}
\begin{gathered}
   \Li_{0}(1/w) = \frac{1/w}{(1-1/w)} = \frac{1}{w-1} 
   = - 1 - \Li_{0}(w),
\\ 
   \Li_{-l}(\frac1w)
   = \left(\frac1w \deriv{(1/w)}\right)^{l} \Li_{0}(\frac1w)
   = - \left(- w\deriv{w}\right)^{l} \Li_{01}(w)
   = (-1)^{l-1} \Li_{-l}(w).
\end{gathered}
\end{equation*}
Next consider $\Li_1(w) = -\log(1-w)$. We need to specify how we take
the branch of $\log$. We set $w = e^{-y}$ with $y > 0$. Then
$\Li_1(e^{-y})$ is defined. We define $\Li_1(e^{-y})$ with $y < 0$ by
analytic continuation. Then
\begin{equation}\label{eq:formula_Li1}
   \Li_1(e^{y}) = \Li_1(e^{-y}) - y - \pi\sqrt{-1}.
\end{equation}

\begin{NB}
Notation: $\arg(y)=0$ if $y>0$. $\arg(-1)=\pi$ ($\arg((-1)^{-1})=-\pi$).
\end{NB}

If $2\pi>y>0$, then by integrating
\begin{equation}
\log(1-e^{-y})=\log y-\frac{1}{2}y+
\frac{B_1}{2 \cdot 2!}y^2-\frac{B_2}{4 \cdot 4!}y^4-\cdots,
\end{equation}
we have
\begin{equation}
\Li_2(e^{-y})=\frac{\pi^2}{6}+(y \log y-y)-\frac{y^2}{4}+
\frac{B_1}{2 \cdot 3 \cdot 2!}y^3-\frac{B_2}{4 \cdot 5 \cdot 4!}y^5-\cdots.
\end{equation}
Note that
\begin{equation}
-\frac{y^2}{4}+
\frac{B_1}{2 \cdot 3 \cdot 2!}y^3-\frac{B_2}{4 \cdot 5 \cdot 4!}y^5-\cdots.
\end{equation}
converges for $|y|<2\pi$.
We define $\Li_2(e^{-y})$, $y<0$ by analytic continuation. 
So 
\begin{equation}\label{eq:Li_2}
\Li_2(e^{-y})=\frac{\pi^2}{6}+(y \log y-y)-\frac{y^2}{4}+
\frac{B_1}{2 \cdot 3 \cdot 2!}y^3-\frac{B_2}{4 \cdot 5 \cdot 4!}y^5-\cdots.
\end{equation}
for $|y|<2\pi$.
\begin{NB}
This definition is the same as
\begin{equation}
\Li_2(w)=-\int_0^w \frac{\log(1-u)}{u}du
\end{equation}
where the branch cut of $\log(1-u)$ is $[1,\infty]$
and $-\pi<\arg(1-u) \leq \pi$.
Thus if $w>1$, then the path is on the lower half plane. 
\end{NB}
Hence
\begin{equation}\label{formula:Li2}
\begin{split}
\Li_2(e^y)+\Li_2(e^{-y})&=\frac{\pi^2}{3}-y\log(-y)+y\log y-\frac{y^2}{2}\\
&=\frac{\pi^2}{3}-y \pi \sqrt{-1}-\frac{y^2}{2}
\end{split}
\end{equation}
for $y>0$.
By integrating \eqref{eq:Li_2}, we have 
\begin{equation}
\Li_3(e^{-y})=\zeta(3)-\frac{\pi^2}{6}y-\frac{y^2}{2}\log y+
\frac{3}{4}y^2+\frac{y^3}{12}
-\left(\frac{B_1}{2 \cdot 3 \cdot 4 \cdot 2!}y^4-
\frac{B_2}{4 \cdot 5  \cdot 6 \cdot 4!}y^6-\cdots \right)
\end{equation}
for $|y|<2\pi$.
Hence we have
\begin{equation}\label{formula:Li3}
\begin{split}
\Li_3(e^y)&=\Li_3(e^{-y})+\frac{\pi^2}{3}y+\frac{1}{2}y^2 \log y-
\frac{1}{2}y^2 \log (-y)-\frac{y^3}{6}\\
&=\Li_3(e^{-y})+\frac{\pi^2}{3}y-
\frac{1}{2}y^2 \pi \sqrt{-1}-\frac{y^3}{6}
\end{split}
\end{equation}
for $y>0$.

\subsection{Limit}\label{subsec:limit}
We have
\begin{equation*}
   \lim_{\beta\to 0} \beta^{k+1} \Li_{-k}(e^{-\bbeta x})
   = x^{-k-1} P_{k}(1) = k! x^{-k-1}
\end{equation*}
for $k\in\Z_{\ge 0}$.

We have
\begin{equation*}
   \Li_{1}(e^{-\bbeta x}) %- \frac{\bbeta x}2
   + \log(\bbeta\Lambda)
   = - \log \left(
     \frac{1 - e^{-\bbeta x}}{\bbeta\Lambda}
%     \frac2{\bbeta\Lambda}\sinh \frac{\bbeta x}2
   \right)
   \xrightarrow[\bbeta\to 0]{}
   - \log\left(\frac{x}\Lambda\right).
\end{equation*}
Then
\begin{equation*}
\begin{split}
   & \int_0^x 
   \Li_{1}(e^{-\bbeta x'}) %- \frac{\bbeta x'}2
   + \log(\bbeta \Lambda) dx'
   =
   -\frac1\bbeta\left(\Li_2(e^{-\bbeta x}) - \frac{\pi^2}6\right)
     %- \frac{\bbeta x^2}4
     + x\log(\bbeta\Lambda)
\\
  & 
  \xrightarrow[\bbeta\to 0]{}
  - \int_0^x \log\left(\frac{x'}\Lambda\right) dx'
  = - x \log\left(\frac{x}\Lambda\right) + x.
\end{split}
\end{equation*}
Furthermore
\begin{equation*}
\begin{split}
   & \int_0^x 
   -\frac1\bbeta\left(\Li_2(e^{-\bbeta x'}) - \frac{\pi^2}6\right)
     %- \frac{\bbeta (x')^2}4
   + x'\log(\bbeta\Lambda) dx'
\\
   & \qquad =
   \frac1{\bbeta^2}\left( \Li_3(e^{-\bbeta x}) - \zeta(3)\right)
   %- \frac{\bbeta x^3}{12}
   + \frac{x^2}2 \log(\bbeta\Lambda)  + \frac{\pi^2 x}{6\bbeta}
\\
  & 
  \xrightarrow[\bbeta\to 0]{}
  \int_0^x \left[- x' \log\left(\frac{x'}\Lambda\right) + x'\right] dx'
  = - \frac12 x^2 \log\left(\frac{x}\Lambda\right) + \frac34 x^2.
\end{split}
\end{equation*}

\end{document}